# SEMILINEAR ELLIPTIC EQUATIONS WITH HARDY POTENTIAL AND SUBCRITICAL SOURCE TERM

**PHUOC-TAI NGUYEN**
DEPARTAMENTO DE MATEMATICAS, PONTIFICIA UNIVERSIDAD CATÓLICA DE CHILE
CASILLA 307, CORREO 2, SANTIAGO, CHILE

Contents

Abstract. Let $\Omega$ be a smooth bounded domain in $\mathbb{R}^N$ and $\delta(x) = \mathrm{dist}\,(x, \partial\Omega)$. Assume $\mu \in \mathbb{R}$, $\nu$ is a nonnegative finite measure on $\partial\Omega$ and $g \in C(\Omega \times \mathbb{R}_+)$. We study positive solutions of

$$(P) \qquad -\Delta u - \frac{\mu}{\delta^2} u = g(x,u) \text{ in } \Omega, \qquad \mathrm{tr}^*(u) = \nu.$$

Here $\mathrm{tr}^*(u)$ denotes the *normalized boundary trace* of $u$ which was recently introduced by M. Marcus and P. T. Nguyen in [16]. We focus on the case $0 < \mu < C_H(\Omega)$ (the Hardy constant for $\Omega$) and provide some qualitative properties of solutions of (P). When $g(x,u) = u^q$ with $q > 1$, we prove that there is a critical value $q^*$ (depending only on $N$, $\mu$) for (P) in the sense that if $1 < q < q^*$ then (P) admits a solution under a smallness assumption on $\nu$, but if $q \geq q^*$ this problem admits no solution with isolated boundary singularity. Existence result is then extended to a more general setting where $g$ is *subcritical* (see (1.27)). We also investigate the case where the $g$ is linear or sublinear and give some existence results for (P).

*Key words:* Hardy potential, Martin kernel, normalized boundary trace, source terms.

*2000 Mathematics Subject Classification:* 35J60, 35J75, 35J10

Email address: nguyenphuoctai.hcmup@gmail.com.

## 1. Introduction

This paper concerns a study of weak solutions of semilinear elliptic equation with Hardy potential and source term

$$-\Delta u - \frac{\mu}{\delta^2}u = g(x,u) \tag{1.1}$$

in a $C^2$ bounded domain $\Omega$, where $\mu \in \mathbb{R}$, $\delta(x) = \mathrm{dist}\,(x, \partial\Omega)$ and $g \in C(\Omega \times \mathbb{R}_+)$, $g(x,0) = 0$.

Henceforth, we will use the notations $L_\mu := \Delta + \frac{\mu}{\delta^2}$ and $(g \circ u)(x) := g(x, u(x))$.

**Definition 1.1.** (i) *A function $u$ is called a (weak) solution (resp. subsolution, supersolution) of* (1.1) *if $u \geq 0$, $u \in L^1_{loc}(\Omega)$, $g \circ u \in L^1_{loc}(\Omega)$ and*

$$-L_\mu u = g \circ u \; (\textit{resp.} -L_\mu u \leq g \circ u, -L_\mu u \geq g \circ u)$$

*in the sense of distributions.*

(ii) *A function $u$ is an $L_\mu$-harmonic function (resp. $L_\mu$-subharmonic, $L_\mu$-superharmonic) if $u \in L^1_{loc}(\Omega)$ and*

$$-L_\mu u = 0 \; (\textit{resp.} -L_\mu u \leq 0, \; -L_\mu u \geq 0)$$

*in the sense of distributions.*

Boundary value problem with measures for (1.1) with $\mu = 0$ and $g \circ u = u^q$, i.e. the problem,

$$-\Delta u = u^q \quad \text{in } \Omega, \qquad u = \nu \quad \text{on } \partial\Omega. \tag{1.2}$$

was first considered by Bidault-Véron and Vivier in [6]. They established estimates involving classical Green and Poisson kernels for $-\Delta$ and applied these estimates to obtain an existence result in the subcritical case, i.e. $1 < q < q_c := \frac{N+1}{N-1}$. Then Bidault-Véron and Yarur [8] reconsidered this type of problem in a more general setting and provided necessary and sufficient condition for the existence of a solution of (1.2). In [11] Chen, Felmer and Véron investigated (1.1) with $\mu = 0$ and $g$ satisfying a subcritical condition. Their approach makes use of Schauder fixed point theorem, essentially based on estimates related to weighted Marcinkiewicz spaces. Recently, Véron et al. [7] provided new criteria for the existence of weak solutions of problem (1.2) and extended those results to the case where $\Delta$ is replaced by $L_\mu$.

When $\mu \neq 0$, the study of (1.1) relies strongly on the investigation of the linear equation

$$-L_\mu u = 0. \tag{1.3}$$

Equation (1.3) with $\mu < 0$, and more generally Schrödinger equations $-\Delta u - V(x)u = 0$ where $V$ is a nonnegative potential, was studied by Ancona [1, 2], Marcus [15], Ancona and Marcus [3] and by Véron and Yarur [21]. The case $\mu > 0$ was considered by Bandle et al. [4, 5], Marcus and Nguyen [16] and by Gkikas and Véron [14] in connection with the optimal constant $C_H(\Omega)$ in Hardy's inequality, namely

$$C_H(\Omega) = \inf_{H^1_0(\Omega)\setminus\{0\}} \frac{\int_\Omega |\nabla u|^2 dx}{\int_\Omega (u/\delta)^2 dx}. \tag{1.4}$$

It is well known (see [10, 17]) that $C_H(\Omega) \in (0, \frac{1}{4}]$ and $C_H(\Omega) = \frac{1}{4}$ when $\Omega$ is convex. Moreover the infimum is achieved if and only if $C_H(\Omega) < 1/4$.

Let $\phi \geq 0$ in $\Omega$ and $p \geq 1$, we denote by $L^p(\Omega; \phi)$ the space of all function $v$ on $\Omega$ satisfying $\int_\Omega |v|\, \phi dx < \infty$. We denote by $\mathfrak{M}(\Omega; \phi)$ the space of Radon measures $\tau$ on $\Omega$ satisfying $\int_\Omega \phi\, d|\tau| < \infty$ and by $\mathfrak{M}^+(\Omega; \phi)$ the nonnegative cone of $\mathfrak{M}(\Omega; \phi)$. When $\phi \equiv 1$, we use the usual notations $\mathfrak{M}(\Omega)$ and $\mathfrak{M}^+(\Omega)$. We also denote by $\mathfrak{M}(\partial\Omega)$ the space of finite measures on $\partial\Omega$ and by $\mathfrak{M}^+(\partial\Omega)$ the nonnegative cone of $\mathfrak{M}(\partial\Omega)$.

Let $G_\mu$ and $K_\mu$ be the Green and the Martin kernels for $-L_\mu$ in $\Omega$ respectively (see [16] for more detail). Denote by $\mathbb{G}_\mu$ and $\mathbb{K}_\mu$ the associated operators defined by

$$\mathbb{G}_\mu[\tau](x) = \int_\Omega G_\mu(x,y)d\tau(y) \quad \forall \tau \in \mathfrak{M}(\Omega). \tag{1.5}$$

$$\mathbb{K}_\mu[\nu](x) = \int_{\partial\Omega} K_\mu(x,z)d\nu(z) \quad \forall \nu \in \mathfrak{M}(\partial\Omega). \tag{1.6}$$

Put

$$\alpha_\pm := \frac{1 \pm \sqrt{1-4\mu}}{2}. \tag{1.7}$$

Let $\lambda_{\mu,1}$ be the first eigenvalue of $-L_\mu$ in $\Omega$ and denote by $\varphi_{\mu,1}$ the corresponding eigenfunction normalized by $\int_\Omega (\varphi_{\mu,1}/\delta)^2 dx = 1$ (see [10]). If $\mu \in (0, C_H(\Omega))$ then $\lambda_{\mu,1} > 0$ and by [12] (see also [19]), there exists a constant $c_1 > 0$ such that

$$c_1^{-1}\delta^{\alpha_+} \le \varphi_{\mu,1} \le c_1\delta^{\alpha_+} \quad \text{in } \Omega. \tag{1.8}$$

For $\beta > 0$, put

$$\Omega_\beta = \{x \in \Omega : \delta(x) < \beta\},\ D_\beta = \{x \in \Omega : \delta(x) > \beta\},\ \Sigma_\beta = \{x \in \Omega : \delta(x) = \beta\}.$$

When dealing with boundary value problem associated to (1.1) with $\mu > 0$ one encounters the following difficulties:

- The first one is due to the fact that every positive $L_\mu$-harmonic function has classical measure boundary trace zero (see [16, Corollary 2.11]). Therefore, classical boundary trace no longer plays a role in describing the boundary behavior of $L_\mu$-harmonic function or solutions of (1.1).
- The second one stems from the invalidity of the classical Keller-Osserman estimate, as well as the lack of a universal upper bound for solutions of (1.1). Moreover, contrast to the case of nonnegative absorption nonlinearity, $\mathbb{K}_\mu[\nu]$ is a subsolution of

$$-L_\mu u = g \circ u \quad \text{in } \Omega, \qquad \operatorname{tr}^*(u) = \nu \tag{1.9}$$

and therefore it is no longer an upper bound for solutions of (1.9).

In order to overcome the first difficulty, we shall employ the notion of *normalized boundary trace* which is defined as follows:

**Definition 1.2.** *A positive function $u$ possesses a* normalized boundary trace *if there exists a measure $\nu \in \mathfrak{M}^+(\partial\Omega)$ such that*

$$\lim_{\beta\to 0} \int_{\Sigma_\beta} \frac{|u - \mathbb{K}_\mu[\nu]|}{\beta^{\alpha_-}} dS = 0. \tag{1.10}$$

*The normalized boundary trace of $u$ is denoted by* $\operatorname{tr}^*(u)$.

In the above definition, we use the notation $dS = d\mathbb{H}_{N-1}$ where $\mathbb{H}_{N-1}$ denotes the Hausdorff measure. This notion is introduced by Marcus and Nguyen [16] in the case $\mu \in (0, C_H(\Omega))$. It is worth mentioning that if $\mu \in (0, C_H(\Omega))$ then $\lambda_{\mu,1} > 0$ and hence $\varphi_{\mu,1}$ is a positive $L_\mu$-superharmonic function. This fact, together with a classical result of Ancona [2], implies the existence of $L_\mu$ harmonic functions and guarantees the validity of Representation theorem (see [16]). Normalized boundary trace turns out to be a more appropriate notion to investigate the problem

$$-L_\mu u + u^q = 0 \quad \text{in } \Omega, \qquad \operatorname{tr}^*(u) = \nu. \tag{1.11}$$

More precisely, when $\mu \in (0, C_H(\Omega))$, they showed that there exists a critical exponent

$$q^* = q^*(N,\mu) := \frac{N+\alpha_+}{N-1-\alpha_-}. \tag{1.12}$$

for (1.11). This means that if $1 < q < q^*$, for every positive finite boundary measure $\nu$, (1.11) admits a unique positive solution while if $q \geq q_*$ there exists no positive solution of (1.11) with $\nu$ being a Dirac measure. Stability result was also discussed in the case $1 < q < q^*$. Problem (1.11) with $u^q$ replaced by a more general nonlinearity $f(u)$ was then investigated by Gkikas and Véron [14] in a slightly different setting. When $f(u) = |u|^{q-1}u$, they provided a necessary and sufficient condition in terms of Besov capacity for solving (1.11) in the supercritical case, i.e. $q \geq q^*$.

Because of the second difficulty, we mainly deal with the minimal solution of (1.9) which possesses several exploitable properties. This solution is constructed due to sub-supersolutions theorem that is established in Section 3. Observe that $\mathbb{K}_\mu[\nu]$ is a subsolution of (1.9). Hence in order to prove the existence of a minimal solution of (1.9), it is sufficient to find a supersolution of (1.9) which dominates $\mathbb{K}_\mu[\nu]$.

Throughout the present paper, we assume that $\mu \in (0, C_H(\Omega))$. We now introduce the definition of solutions of (1.9).

**Definition 1.3.** (i) *A nonnegative function $u$ is called a* (weak) solution of (1.9) *if $u$ is a solution of* (1.1) *and has normalized boundary trace $\nu$.*

(ii) *Let us define the space of* admissible test function *as follows:*

$$X(\Omega) = \{\zeta \in C^2(\Omega) : \delta^{\alpha_-} L_\mu \zeta \in L^\infty(\Omega),\ \delta^{-\alpha_+}\zeta \in L^\infty(\Omega)\}.$$

*A function $\zeta \in X(\Omega)$ is called an* admissible test function *for* (1.9).

Notice that $\varphi_{\mu,1} \in X(\Omega)$. More properties of $X(\Omega)$ can be found in [16, Section 2.4]. Using this space, we establish integral formulation for weak solutions of (1.9). This is stated in the following result.

**Theorem A.** *Let $\nu \in \mathfrak{M}^+(\partial\Omega)$. The following statements are equivalent:*

(i) *$u$ is a solution of* (1.9).

(ii) *$g \circ u \in L^1(\Omega; \delta^{\alpha_+})$ and*

$$u = \mathbb{G}_\mu[g \circ u] + \mathbb{K}_\mu[\nu]. \tag{1.13}$$

(iii) *$u \in L^1(\Omega; \delta^{-\alpha_-})$, $g \circ u \in L^1(\Omega; \delta^{\alpha_+})$ and*

$$-\int_\Omega u L_\mu \zeta dx = \int_\Omega (g \circ u)\zeta dx - \int_\Omega \mathbb{K}_\mu[\nu] L_\mu \zeta dx \quad \forall \zeta \in X(\Omega). \tag{1.14}$$

Under some additional assumptions on $g$, we obtain existence result for (1.9).

**Theorem B.** *Let $g(x,r)$ be a nondecreasing continuous function with respect to $r$ for every $x \in \Omega$ and $\nu \in \mathfrak{M}^+(\partial\Omega)$ with $\|\nu\|_{\mathfrak{M}(\partial\Omega)} = 1$. Assume that there exist numbers $c_2 > 0$, $c_3 > 0$, $0 \leq r_1 < r_2 \leq \infty$ and a positive function $\ell$ such that*

$$g(x, rs) \leq \ell(r) g(x,s) \quad \forall s \geq 0, r > 0, x \in \Omega. \tag{1.15}$$

$$\ell(1 + c_2 c_3 r^{-1} \ell(r)) \leq c_2 \quad \forall r \in (r_1, r_2), \tag{1.16}$$

$$\mathbb{G}_\mu[g \circ (\mathbb{K}_\mu[\nu])] \leq c_3 \mathbb{K}_\mu[\nu] \quad a.e.\ in\ \Omega. \tag{1.17}$$

1. EXISTENCE. *For any $\varrho \in (r_1, r_2)$ the problem*

$$-L_\mu u = g \circ u \quad \text{in } \Omega, \qquad \operatorname{tr}^*(u) = \varrho\nu \tag{1.18}$$

*admits a minimal solution $\underline{u}_{\varrho\nu}$ in the sense that if $v$ is a solution of (1.18) then $\underline{u}_{\varrho\nu} \leq v$.*
2. ESTIMAES. *This solution satisfies*

$$\varrho\mathbb{K}_\mu[\nu] \leq \underline{u}_{\varrho\nu} \leq c_4\varrho\mathbb{K}_\mu[\nu] \tag{1.19}$$

*where $c_4 = c_4(c_2, c_3, \ell, \varrho)$.*
3. NONTANGENTIAL CONVERGENCE. *For $\nu$-a.e. $z \in \partial\Omega$, there holds*

$$\lim_{x\to z} \frac{\underline{u}_{\varrho\nu}(x)}{\mathbb{K}_\mu[\nu](x)} = \varrho \quad \text{non tangentially.} \tag{1.20}$$

**Remark.** As a model, we can take $g(x, u) = \ln^\beta(u + 1)u^q$ with small $\beta \geq 0$ and $q > 0$.

In the next results, we focus on the pure power case, namely the problem

$$(D_\nu) \qquad -L_\mu u = u^q \quad \text{in } \Omega, \qquad \operatorname{tr}^*(u) = \nu$$

where $q > 0$ and $\nu \in \mathfrak{M}^+(\partial\Omega)$. We shall establish some estimates related Green and Martin operators and a necessary condition for the existence of solutions of $(D_\nu)$ in the case $q > 1$.

**Proposition C.** *Let $q > 0$ and $\nu \in \mathfrak{M}^+(\partial\Omega)$. Then there exists a positive constant $c_5 = c_5(N, \mu, q, \Omega)$ such that*

$$\mathbb{G}_\mu[\mathbb{K}_\mu[\nu]^q] \leq c_5 \|\nu\|^{q-1}_{\mathfrak{M}(\partial\Omega)} \mathbb{K}_\mu[\nu]. \tag{1.21}$$

*Furthermore, if $q > 1$ and problem $(D_\nu)$ admits a solution then there holds*

$$\mathbb{G}_\mu[\mathbb{K}_\mu[\nu]^q] \leq \frac{1}{q-1}\mathbb{K}_\mu[\nu]. \tag{1.22}$$

**Remark.** It is worth mentioning that (1.21) with $q > 1$ and (1.22) with an inexplicit multiplier were proved in [7, Theorem 4.1]. In this paper we employ the method in [8] to prove (1.21) for $q > 0$ and apply the idea in [9] to point out that the multiplier can be explicitly chosen as $\frac{1}{q-1}$.

The next results reveal that $q^*$ is the *critical exponent* for $(D_\nu)$. More precisely, in the subcritical case, namely $1 < q < q^*$, $(D_\nu)$ admits a solution under smallness assumption on the boundary datum while in the supercritical case, i.e. $q \geq q^*$, this problem possesses no solution with isolated boundary singularity.

For $z \in \partial\Omega$, we denote by $\delta_z$ the Dirac measure concentrated at $z$. Existence and nonexistence results when $0 < q < q^*$ , $q \neq 1$ is given as follows.

**Theorem D.** *Let $q \in (0, q^*)$, $q \neq 1$ and $\nu \in \mathfrak{M}^+(\partial\Omega)$ with $\|\nu\|_{\mathfrak{M}(\partial\Omega)} = 1$. For $\varrho > 0$, consider the problem*

$$(D_{\varrho\nu}) \qquad -L_\mu u = u^q \quad \text{in } \Omega, \qquad \operatorname{tr}^*(u) = \varrho\nu.$$

1. CASE $q > 1$. *There exists a threshold value $\varrho^* \in \mathbb{R}_+$ for $(D_{\varrho\nu})$ such that*

(i) *If $\varrho \in (0, \varrho^*]$ then problem $(D_{\varrho\nu})$ admits a minimal solution $\underline{u}_{\varrho\nu}$. Moreover, if $\varrho \in (0, \varrho^*)$, $\underline{u}_{\varrho\nu}$ satisfies (1.19) and (1.20).*

*In addition, if $\{\varrho_n\}$ be a nondecreasing sequence converging to $\varrho^*$ then the sequence $\{\underline{u}_{\varrho_n\nu}\}$ converges to $\underline{u}_{\varrho^*\nu}$ in $L^1(\Omega; \delta^{-\alpha_-})$ and in $L^q(\Omega; \delta^{\alpha_+})$.*

(ii) *If $\varrho > \varrho^*$ then there exists no solution of $(D_{\varrho\nu})$.*

2. Case $q \in (0,1)$. *For every $\varrho > 0$ problem $(D_{\varrho\nu})$ admits a minimal solution $\underline{u}_{\varrho\nu}$ which satisfies satisfies* (1.19) *and* (1.20)*. Moreover,* $\lim_{\varrho\to 0}\underline{u}_{\varrho\nu} = 0$ *and* $\lim_{\varrho\to\infty}\underline{u}_{\varrho\nu} = \infty$.

*In any case, if $\nu = \delta_z$ with $z \in \partial\Omega$ then there holds*

$$\lim_{x\to z}\frac{\underline{u}_{\varrho\delta_z}(x)}{K_\mu(x,z)} = \varrho. \tag{1.23}$$

**Remark.** It is worth noticing that in absorption case (1.11), if $1 < q < q^*$, there are two types of solution with isolated boundary singularity: weakly singular solutions $u_{\varrho,z}$ (the solution of (1.11) with $\nu = \varrho\delta_z$) and strongly singular solution $u_{\infty,z}$. Actually, $u_{\infty,z}$ is the limit of the sequence $u_{\varrho,z}$ as $\varrho \to \infty$. This limiting process can not be executed in the source case since $(D_{\varrho\delta_z})$ admits no solution if $\varrho > \varrho^*$ due to Theorem D.

We next give a stability result.

**Theorem E.** *Let $q \in (0,q^*)$, $q \neq 1$ and $\{\nu_n\}$ is a sequence of measures in $\mathfrak{M}^+(\partial\Omega)$ which converges weakly to $\nu \in \mathfrak{M}^+(\partial\Omega)$. If $q > 1$, assume, in addition, that*

$$\sup_n \|\nu_n\|_{\mathfrak{M}(\partial\Omega)} \leq \varrho^*. \tag{1.24}$$

*For each $n$, let $u_{\nu_n}$ be a solution of $(D_{\nu_n})$. Then, up to a subsequence, $\{u_{\nu_n}\}$ converges to a solution $u_\nu$ of $(D_\nu)$ in $L^1(\Omega;\delta^{-\alpha_-})$ and in $L^q(\Omega;\delta^{\alpha_+})$.*

Existence and stability result in the case $q = 1$ is stated in the following theorem in which $\lambda_{\mu,1}$ is the first eigenvalue of $-L_\mu$ in $\Omega$.

**Theorem F.** *Let $\nu \in \mathfrak{M}^+(\partial\Omega)$. For $\kappa > 0$, consider the problem*

$$(E_\nu^\kappa) \qquad -L_\mu u = \kappa u \quad \text{in } \Omega; \qquad \operatorname{tr}^*(u) = \nu.$$

*There exists a number $\kappa^* \in (0,\lambda_{\mu,1}]$ such that*

(i) *If $\kappa \in (0,\kappa^*)$ then problem $(E_\nu^\kappa)$ admits a minimal solution $\underline{u}_{\kappa,\nu}$. Moreover, $\underline{u}_{\kappa,\nu}$ satisfies* (1.20).

*Assume $\{\nu_n\}$ is a sequence of measures in $\mathfrak{M}^+(\partial\Omega)$ which converges weakly to $\nu \in \mathfrak{M}^+(\partial\Omega)$ and for each $n$ denote by $u_{\kappa,\nu_n}$ a solution of $(E_{\nu_n}^\kappa)$. Then, up to a subsequence, $\{u_{\kappa,\nu_n}\}$ converges to a solution $u_{\kappa,\nu}$ of $(E_\nu^\kappa)$ in $L^1(\Omega;\delta^{-\alpha_-})$.*

(ii) *If $\kappa > \kappa^*$ then $(E_\nu^\kappa)$ admits no solution.*

*Furthermore, problem $(E_\nu^{\lambda_{\mu,1}})$ admits no solution.*

Note that the assumption that $g(x,r)$ is nondecreasing with respect to $r$ is crucial to obtain the existence in Theorems D and F. A natural question arises: "Does the existence results still hold if *the monotonicity condition is dropped*?". Positive answer to this question is given in the next two theorems where a more general weighted source term $g(x,r)$ is involved. More precisely, we consider the case $(g\circ u)(x) = \delta(x)^\gamma \tilde{g}(u(x))$ with $\gamma > -1-\alpha_+$ and $\tilde{g}: \mathbb{R}_+ \to \mathbb{R}_+$ being continuous. In this framework, the critical value is defined as follows:

$$q_\gamma^* = q_\gamma^*(N,\mu,\gamma) := \frac{N+\alpha_+ +\gamma}{N-1-\alpha_-}. \tag{1.25}$$

Clearly $q_0^* = q^*$.

Theorem G gives existence result for the problem

$$-L_\mu u = \delta^\gamma \tilde{g}(u) \quad \text{in } \Omega, \qquad \operatorname{tr}^*(u) = \varrho\nu. \tag{1.26}$$

**Theorem G.** *Let $\nu \in \mathfrak{M}^+(\partial\Omega)$ such that $\|\nu\|_{\mathfrak{M}(\partial\Omega)} = 1$. Assume that*

$$\Lambda_0 := \int_1^\infty s^{-1-q^*_\gamma} \tilde{g}(s) ds < +\infty, \tag{1.27}$$

$$\tilde{g}(s) \le \Lambda_1 s^{q_1} + \theta \quad \forall s \in [0,1] \quad \textit{for some } q_1 > 1, \Lambda_1 > 0, \theta > 0. \tag{1.28}$$

*There exist $\theta_0 > 0$ and $\varrho_0 > 0$ depending on $N$, $\mu$, $\gamma$, $\Lambda_0$, $\Lambda_1$ and $q_1$ such that for every $\theta \in (0, \theta_0)$ and $\varrho \in (0, \varrho_0)$ problem (1.26) admits a nonnegative solution.*

**Remark.** We say that $g$ is *subcritical* if $\tilde{g}$ satisfies (1.27).

The case where $g$ is *linear* or *sublinear* is treated in the following theorem.

**Theorem H.** *Let $\nu \in \mathfrak{M}^+(\partial\Omega)$ such that $\|\nu\|_{\mathfrak{M}(\partial\Omega)} = 1$. Assume that*

$$\tilde{g}(s) \le \Lambda_2 s^{q_2} + \theta \quad \forall s \ge 0 \tag{1.29}$$

*for some $q_2 \in (0, 1]$, $\Lambda_2 > 0$ and $\theta > 0$.*

*In (1.29), if $q_2 = 1$ we assume in addition that $\Lambda_2$ is small enough. Then for any $\varrho > 0$, (1.26) admits a nonnegative solution.*

Note that in Theorem H, when $q_2 < 1$, smallness assumption on $\theta$ is not required.

When $\tilde{g}$ does not satisfy (1.27), there is no solution with an isolated boundary singularity. This is stated in the following Theorem where we assume that $\tilde{g}$ is nondecreasing.

**Theorem I.** *Assume $\tilde{g}$ is a nondecreasing function such that*

$$\int_a^\infty s^{-1-q^*_\gamma} \tilde{g}(s) ds = \infty \quad \textit{for some } a > 0. \tag{1.30}$$

*Then for every $\varrho > 0$ and $z \in \partial\Omega$ there exists no positive solution of*

$$-L_\mu u = \delta^\gamma \tilde{g}(u) \quad \textit{in } \Omega, \qquad \operatorname{tr}^*(u) = \varrho \delta_z \tag{1.31}$$

**Remark.** If $g(u) = u^q$ then (1.30) is satisfied if and only if $q$ belongs to supercritical range, i.e. $q \ge q^*_\gamma$. We notice that Theorem I was obtained in [7] for the case $\gamma = 0$ and $g(u) = u^q$. Interesting existence results for $(D_{\rho\nu})$ in the supercritical case are also provided in [7].

The plan of the paper is as follows. In Section 2 we give some results concerning Green and Martin kernels and boundary value problem for linear equations with Hardy potential. Theorems A and B are proved in Section 3. It is noteworthy that main ingredients in proving Theorem A is a generalization of Herglotz-Doob to $L_\mu$-superharmonic functions and theory of Schrödinger linear equations. Theorem B is established using a sub-supersolutions theorem. The proof of Proposition C and Theorems D-F are provided in Section 4. Finally, in Section 5 we present the proof of the existence result for a more general class of source terms (Theorems G and H) and demonstrate the nonexistence result in the supercritical case (Theorem I).

**Acknowledgements** The author is grateful to Q. H. Nguyen for his useful comments.

## 2. Preliminaries

Throughout this paper we assume that $0 < \mu < C_H(\Omega)$.

2.1. **Weak $L^p$ spaces.** Denote $L^p_w(\Omega;\tau)$, $1 \le p < \infty$, $\tau \in \mathfrak{M}^+(\Omega)$, the weak $L^p$ space (or Marcinkiewicz space) defined as follows: a measureable function $f$ in $\Omega$ belongs to this space if there exists a constant $c$ such that

$$\lambda_f(a;\tau) := \tau(\{x \in \Omega : |f(x)| > a\}) \le c a^{-p}, \quad \forall a > 0. \tag{2.1}$$

The function $\lambda_f$ is called the distribution function of $f$ (relative to $\tau$). For $p \ge 1$, denote

$$L^p_w(\Omega;\tau) = \{f \text{ Borel measurable} : \sup_{a>0} a^p \lambda_f(a;\tau) < \infty\},$$

$$\|f\|^*_{L^p_w(\Omega;\tau)} = (\sup_{a>0} a^p \lambda_f(a;\tau))^{\frac{1}{p}}. \tag{2.2}$$

The $\|.\|_{L^p_w(\Omega;\tau)}$ is not a norm, but for $p > 1$, it is equivalent to the norm

$$\|f\|_{L^p_w(\Omega;\tau)} = \sup\left\{ \frac{\int_\omega |f| d\tau}{\tau(\omega)^{1/p'}} : \omega \subset \Omega, \omega \text{ measurable}, \ 0 < \tau(\omega) < \infty \right\}. \tag{2.3}$$

More precisely,

$$\|f\|^*_{L^p_w(\Omega;\tau)} \le \|f\|_{L^p_w(\Omega;\tau)} \le \frac{p}{p-1} \|f\|^*_{L^p_w(\Omega;\tau)}. \tag{2.4}$$

When $\tau = \delta^\alpha dx$, for simplicity, we use the notation $L^p_w(\Omega;\delta^\alpha)$. Notice that, for every $\alpha > -1$,

$$L^p_w(\Omega;\delta^\alpha dx) \subset L^r(\Omega;\delta^\alpha) \quad \forall r \in [1,p).$$

From (2.2) and (2.4), one can derive the following estimate which is useful in the sequel. If $u \in L^p_w(\Omega;\delta^\alpha)$ $(\alpha > -1)$ then

$$\int_{|\{u| \ge s\}} \delta^\alpha dx \le s^{-p} \|u\|^p_{L^p_w(\Omega;\delta^\alpha)}. \tag{2.5}$$

2.2. **Green and Martin kernels.** Let $G_\mu$ be the Green kernel for the operator $-L_\mu$ in $\Omega \times \Omega$ and denote by $\mathbb{G}_\mu$ the associated operator defined by (1.5). It was shown in [16] that for every $\tau \in \mathfrak{M}(\Omega;\delta^{\alpha_+})$, $|\mathbb{G}_\mu[\tau]| < \infty$ a.e. in $\Omega$. Denote by $K_\mu$ the Martin kernel for $-L_\mu$ in $\Omega$ and by $\mathbb{K}_\mu$ the Martin operator defined by (1.6).

In what follows the notation $f \sim g$ means: there is a positive constant $c$ such that $c^{-1} f < g < cf$ in the domain of the two functions.

By [13, Theorem 4.11] and [16] (see also [14]),

$$G_\mu(x,y) \sim \min\left\{ |x-y|^{2-N}, \delta(x)^{\alpha_+} \delta(y)^{\alpha_+} |x-y|^{2\alpha_- - N} \right\}, \ \forall x,y \in \Omega, \ x \ne y. \tag{2.6}$$

$$K_\mu(x,z) \sim \delta(x)^{\alpha_+} |x-z|^{2\alpha_- - N}, \ \forall x \in \Omega, z \in \partial\Omega. \tag{2.7}$$

By combining (2.6), (2.7) and the estimates of [20, Lemma 2.3.2], we obtain the following.

**Proposition 2.1.** *Let $\varsigma \in [0,1]$. There exist constants $c_i = c_i(N,\mu,\beta,\Omega)$ $(i = 6,7,8)$ such that,*

$$\|\mathbb{G}_\mu[\tau]\|_{L_w^{\frac{N+\beta}{N-2+2\varsigma\alpha_+}}(\Omega,\delta^{\beta-\varsigma\alpha_+})} \le c_6 \|\tau\|_{\mathfrak{M}(\Omega;\delta^{\varsigma\alpha_+})}, \ \forall \tau \in \mathfrak{M}(\Omega;\delta^{\varsigma\alpha_+}), \tag{2.8}$$

*with* $\beta > \max\{-1, -2 + 2\varsigma\alpha_+\}$,

$$\|\mathbb{G}_\mu[\tau]\|_{L_w^{\frac{N+\beta}{N-1-\alpha_-}}(\Omega,\delta^\beta)} \le c_7 \|\tau\|_{\mathfrak{M}(\Omega;\delta^{\alpha_+})}, \quad \forall \tau \in \mathfrak{M}(\Omega;\delta^{\alpha_+}), \ \beta > -1. \tag{2.9}$$

$$\|\mathbb{K}_\mu[\nu]\|_{L_w^{\frac{N+\beta}{N-1-\alpha_-}}(\Omega,\delta^\beta)} \le c_8 \|\nu\|_{\mathfrak{M}(\partial\Omega)}, \quad \forall \nu \in \mathfrak{M}(\partial\Omega), \ \beta > -1. \tag{2.10}$$

**Proof.** By (2.6), for every $\varsigma \in [0,1]$ there exists a constant $c_6'$ such that

$$G_\mu(x,y) \le c_6' \delta(x)^{\varsigma\alpha_+}\delta(y)^{\varsigma\alpha_+}|x-y|^{2-N-2\varsigma\alpha_+} \quad \forall x \neq y.$$

By proceeding as in the proof of [20, Lemma 2.3.3], we obtain (2.8).

We next prove (2.9). If $\delta(y) \le 2|x-y|$ then

$$G_\mu(x,y) \le c_7' \delta(x)^{\alpha_+}\delta(y)^{\alpha_+}|x-y|^{2-N} \le c_7' 2^{\alpha_+}\delta(x)^{\alpha_+}|x-y|^{1+\alpha_- -N} \quad \forall x \neq y. \tag{2.11}$$

If $\delta(y) > 2|x-y|$ then from the inequality $|\delta(x)-\delta(y)| \le |x-y|$ we deduce that $\delta(x) \ge \delta(y) - |x-y| > |x-y|$. Therefore

$$G_\mu(x,y) \le c_7''|x-y|^{2-N} \le c_7''\delta(x)^{\alpha_+}|x-y|^{1+\alpha_- -N}. \tag{2.12}$$

Combining (2.11) and (2.12) yields

$$G_\mu(x,y) \le c_8'\delta(x)^{\alpha_+}|x-y|^{1+\alpha_- -N} \quad \forall x \neq y \tag{2.13}$$

where $c_8' = \max\{c_7'', c_7' 2^{\alpha_+}\}$. Using similar argument as in the proof of [20, Lemma 2.3.3], we deduce (2.9).

Estimate (2.10) was already obtained in [16, Proposition 2.8].

2.3. **Some result on linear equations.** In this subsection, we recall some results concerning boundary value problem for non-homogeneous linear equation.

$$-L_\mu u = \tau. \tag{2.14}$$

**Definition 2.2.** (i) *A function $u$ is a solution of (2.14) if $u \in L^1_{loc}(\Omega)$ and (2.14) is understood in the sense of distributions.*

(ii) *Let $\tau \in \mathfrak{M}(\Omega;\delta^{\alpha_+})$ and $\nu \in \mathfrak{M}^+(\partial\Omega)$. A function $u$ is a solution of*

$$-L_\mu u = \tau \quad \text{in } \Omega, \qquad \operatorname{tr}^*(u) = \nu, \tag{2.15}$$

*if $u$ is a solution of (2.14) and $u$ admits normalized boundary trace $\nu$.*

**Definition 2.3.** *A nonnegative $L_\mu$-superharmonic function is called an $L_\mu$-potential if its largest $L_\mu$-harmonic minorant is zero.*

The following results, which can be found in [16, Theorem I], is crucial in proving Theorem A.

**Proposition 2.4.** (i) *If $\tau = 0$ then problem (2.15) has a unique solution, $u = \mathbb{K}_\mu[\nu]$. If $u$ is a non-negative $L_\mu$-harmonic function and $\operatorname{tr}^*(u) = 0$ then $u = 0$.*

(ii) *If $\tau \in \mathfrak{M}^+(\Omega;\delta^{\alpha_+})$ then $\mathbb{G}_\mu[\tau]$ has normalized trace zero. Thus $\mathbb{G}_\mu[\tau]$ is a solution of (2.15) with $\nu = 0$.*

(iii) *Let $u$ be a positive $L_\mu$-subharmonic function. If $u$ is dominated by an $L_\mu$-superharmonic function then $L_\mu u \in \mathfrak{M}^+(\Omega;\delta^{\alpha_+})$ and $u$ has a normalized boundary trace. In this case $\operatorname{tr}^*(u) = 0$ if and only if $u \equiv 0$.*

(iv) *Let $u$ be a positive $L_\mu$-superharmonic function. Then there exist $\nu \in \mathfrak{M}^+(\partial\Omega)$ and $\tau \in \mathfrak{M}^+(\Omega;\delta^{\alpha_+})$ such that*

$$u = \mathbb{G}_\mu[\tau] + \mathbb{K}_\mu[\nu]. \tag{2.16}$$

*In particular, $u$ is an $L_\mu$-potential if and only if $\operatorname{tr}^*(u) = 0$.*

(v) *For every $\nu \in \mathfrak{M}^+(\partial\Omega)$ and $\tau \in \mathfrak{M}^+(\Omega;\delta^{\alpha_+})$, problem (2.15) has a unique positive solution. The solution is given by (2.16). Moreover, there exists a positive constant $c_9 = c_9(N,\mu,\Omega)$ such that*

$$\|u\|_{L^1(\Omega;\delta^{-\alpha_-})} \le c_9(\|\tau\|_{\mathfrak{M}(\Omega;\delta^{\alpha_+})} + \|\nu\|_{\mathfrak{M}(\partial\Omega)}). \tag{2.17}$$

(vi) *$u$ is a solution of of* (2.15) *if and only if $u \in L^1(\Omega;\delta^{-\alpha_-})$ and*

$$-\int_\Omega u L_\mu \zeta dx = \int_\Omega \zeta d\tau - \int_\Omega \mathbb{K}_\mu[\nu] L_\mu \zeta dx \quad \forall \zeta \in X(\Omega). \tag{2.18}$$

For easy reference, we present a potential theoretic result which serves to prove Theorem B.

**Theorem 2.5.** *Let $w_1$ be a positive $L_\mu$-potential and $w_2$ be a positive $L_\mu$-harmonic function with $\nu = \operatorname{tr}^*(w_2)$. Assume that $\frac{w_1}{w_2}$ satisfies the local Harnack inequality. Then for $\nu$-a.e. $z \in \partial\Omega$,*

$$\lim_{x\to z} \frac{w_1(x)}{w_2(x)} = 0 \quad \textit{non-tangentially.}$$

This Proposition Can be obtained by combining the Fatou convergence theorem [1, Theorem 1.8] and the fact that if a function satisfies the Harnack inequality, fine convergence at the boundary (in the sense of [1]) implies non-tangential convergence (for more details, see [3]).

## 3. Nonlinear equations with source term

In this section, we deal with nonlinear equations involving source term

$$-L_\mu u = g \circ u \tag{3.1}$$

in $\Omega$ where $0 < \mu < C_H(\Omega)$ and $g : \Omega \times \mathbb{R}_+ \to \mathbb{R}_+$ is continuous.

### 3.1. Equivalent formulation.

For $z \in \partial\Omega$, denote by $\mathbf{n}_z$ the outward unit normal vector to $\partial\Omega$ at $z$. We recall below a geometric property of $C^2$ domains (see [20]).

**Proposition 3.1.** *There exists $\beta_0 > 0$ such that for every point $x \in \overline{\Omega}_{\beta_0}$, there exists a unique point $\sigma_x \in \partial\Omega$ such that $x = \sigma_x - \delta(x)\mathbf{n}_{\sigma_x}$. The mappings $x \mapsto \delta(x)$ and $x \mapsto \sigma_x$ belong to $C^2(\overline{\Omega}_{\beta_0})$ and $C^1(\overline{\Omega}_{\beta_0})$ respectively. Moreover, $\lim_{x\to\sigma(x)} \nabla\delta(x) = -\mathbf{n}_{\sigma_x}$.*

For $D \Subset \Omega$, let $G_\mu^D$ and $K_\mu^D$ be the Green and Poisson kernels of $-L_\mu$ in $D$ respectively. Denote by $\mathbb{G}_\mu^D$ and $\mathbb{K}_\mu^D$ the corresponding Green and Poisson operators in $D$.

We prove below main properties of solutions of (1.9).

**Proof of Theorem A.**

**(i) $\Longrightarrow$ (ii).** Assume $u$ is a positive solution of (1.9). Put $\tau = g \circ u$ and denote $\tau_\beta = \tau|_{D_\beta}$ and $\lambda_\beta = u|_{\Sigma_\beta}$ for $\beta \in (0, \beta_0)$. Consider the boundary value problem

$$-L_\mu v = \tau_\beta \quad \text{in } D_\beta, \qquad v = \lambda_\beta \quad \text{on } \Sigma_\beta.$$

This problem admits a unique solution $v_\beta$ (the uniqueness is derived from [5, Lemma 2.1] since $\mu < C_H(\Omega)$). Therefore $v_\beta = u|_{D_\beta}$. We have

$$u|_{D_\beta} = v_\beta = \mathbb{G}_\mu^{D_\beta}[\tau_\beta] + \mathbb{K}_\mu^{D_\beta}[\lambda_\beta].$$

It follows that

$$\int_{D_\beta} G_\mu^{D_\beta}(\cdot, y)(g \circ u)(y)\, dy = \mathbb{G}_\mu^{D_\beta}[\tau_\beta] \le u|_{D_\beta}.$$

Letting $\beta \to 0$, we get

$$\int_\Omega G_\mu(\cdot, y)(g \circ u)(y)\, dy \le u. \tag{3.2}$$

Fix a point $x_0 \in \Omega$ such that $u(x_0) < \infty$. Keeping in mind that $G(x_0, y) > c_{x_0}\delta(y)^{\alpha_+}$ for every $y \in \Omega$, we deduce from (3.2) that $g \circ u \in L^1(\Omega; \delta^{\alpha_+})$. Thanks to Proposition 2.4 (v), we obtain (1.13).

**(ii) $\Longrightarrow$ (i).** Assume $u$ is a function such that $g \circ u \in L^1(\Omega; \delta^{\alpha_+})$ and (1.13) holds. By Proposition 2.4 (i) $-L_\mu \mathbb{K}_\mu[\nu] = 0$, which implies that $u$ is a solution of (3.1). On the other hand, since $g \circ u \in L^1(\Omega; \delta^{\alpha_+})$, we deduce from Proposition 2.4 (ii) that $\operatorname{tr}^*(\mathbb{G}_\mu[g \circ u]) = 0$. Consequently, $\operatorname{tr}^*(u) = \operatorname{tr}^*(\mathbb{K}_\mu[\nu]) = \nu$.

**(i) $\Longrightarrow$ (iii).** Assume $u$ is a positive solution of (1.9). From the implication (i) $\Longrightarrow$ (ii), we deduce that $u \in L^1(\Omega; \delta^{-\alpha_-})$ and $g \circ u \in L^1(\Omega; \delta^{\alpha_+})$. Hence, by Proposition 2.4 (vi), $u$ satisfies (1.14).

**(iii) $\Longrightarrow$ (i).** This implication follows straightfoward from Proposition 2.4 (vi). □

3.2. **Nondecreasing source.** We start with an existence result for (3.1) in presence of sub and super solutions.

**Theorem 3.2.** *Let $g \in C(\Omega \times \mathbb{R}_+)$, $g(x, r)$ be nondecreasing with respect to $r$ for any $x \in \Omega$. Assume that there exist a subsolution $V_1$ and a supersolution $V_2$ of (3.1) such that $0 \le V_1 \le V_2$ in $\Omega$. Then there exists a solution $u$ of (3.1) which satisfies $V_1 \le u \le V_2$ in $\Omega$.*

*Moreover, if $V_1 = \mathbb{K}_\mu[\nu]$ for some $\nu \in \mathfrak{M}^+(\partial\Omega)$ and $g \circ V_2 \in L^1(\Omega; \delta^{\alpha_+})$ then $u$ is the minimal solution of (1.9) in the sense that $u \le v$ for every solution $v$ of (1.9).*

**Lemma 3.3.** *Let $D \Subset \Omega$, $f \in L^1(D)$, $f \ge 0$ and $\eta \in L^1(\partial D)$, $\eta \ge 0$. Then there exists a unique solution of*

$$-L_\mu u = f \quad \text{in } D, \qquad u = \eta \quad \text{on } \partial D. \tag{3.3}$$

**Proof.** We start with the case $f \in L^\infty(D)$ and $\eta = 0$. Let us consider the functional

$$\mathcal{J}(v) = \frac{1}{2}\int_D \left(|\nabla v|^2 - \frac{\mu}{\delta^2} v^2\right) dx - \int_D f v dx$$

over the sapce $H_0^1(D)$. It can be verified that $\mathcal{J}$ is a convex and semicontinuous on $H_0^1(D)$. Furthermore, since $\mu < C_H(\Omega)$, it follows that the functional $\mathcal{J}$ is coercive. Therefore the problem $\min_{H_0^1(D)} \mathcal{J}(v)$ admits a solution $v \in H_0^1(D)$. The minimizer $v$ is the unique weak solution of (3.3).

If $f \in L^1_+(D)$ then we can approximate it by an increasing sequence $\{f_m\} \subset L^\infty(D)$. Let $v_m$ be the solution of (3.3) with $\eta = 0$ and $f$ replaced by $f_m$. By comparison principle [5, Lemma 2.1], $\{v_m\}$ increases and therefore $v := \lim_{m\to\infty} v_m$ is a solution of (3.3) with $\eta = 0$.

We next consider the case $\eta \in L^1(\partial D)$. Let $v$ be a solution of (3.3) with $\eta = 0$ then $u = v + \mathbb{P}_\mu^D[\eta]$ is a solution of (3.3). The uniqueness follows from the comparison principle. □

**Proof of Theorem 3.2.** Put $u_0 = V_1$ and $\eta_\beta = V_1|_{\Sigma_\beta}$ for $\beta \in (0, \beta_0)$. For $n \ge 1$, consider the problem

$$-L_\mu u = g \circ u_{n-1} \quad \text{in } D_\beta, \qquad u = \eta_\beta \quad \text{on } \partial D_\beta. \tag{3.4}$$

For each $n \ge 1$, by Lemma 3.3 there exists a unique solution $u_{\beta,n}$ of (3.4). Moreover, since $g(x, r)$ is nondecreasing with respect to $r$ for every $x \in \Omega$, by applying comparison principle, we deduce that

$$V_1 \le u_{\beta,n} \le u_{\beta,n+1} \le V_2$$

in $D_\beta$. Therefore $u_\beta := \lim_{n\to\infty} u_{\beta,n}$ is a solution of (3.1) in $D_\beta$ which satisfies $V_1 \le u_\beta \le V_2$ in $D_\beta$. Moreover,

$$u_\beta = \mathbb{G}_\mu^{D_\beta}[g \circ u_\beta] + \mathbb{P}_\mu^{D_\beta}[\eta_\beta]. \tag{3.5}$$

For $0 < \beta' < \beta < \beta_0$, by the comparison principle, $u_{\beta,1} \le u_{\beta',1}$ in $D_\beta$. By the monotonicity assumption on $g$, it follows that $u_{\beta,n} \le u_{\beta',n}$ in $D_\beta$ for every $n > 1$. Therefore $V_1 \le u_\beta \le u_{\beta'} \le V_2$ in $D_\beta$ and hence $u := \lim_{\beta\downarrow 0} u_\beta$ is a solution of (3.1) in $\Omega$ satisfying $V_1 \le u \le V_2$.

In the case $V_1 = \mathbb{K}_\mu[\nu]$, formulation (3.5) becomes

$$u_\beta = \mathbb{G}_\mu^{D_\beta}[g \circ u_\beta] + \mathbb{K}_\mu[\nu]. \tag{3.6}$$

Since $0 \le g \circ u_\beta \le g \circ V_2 \in L^1(\Omega; \delta^{\alpha_+})$, it follows that

$$\lim_{\beta\downarrow 0} \mathbb{G}_\mu^{D_\beta}[g \circ u_\beta] = \mathbb{G}_\mu[g \circ u].$$

Letting $\beta \downarrow 0$ in (3.6), we infer that $u$ satisfies (1.13), namely $u$ is a solution of (1.9). Notice that in the above argument, $u$ is independent of $V_2$. Hence, if $v$ is a solution of (1.9) then $v \ge \mathbb{K}_\mu[\nu]$ and $g \circ v \in L^1(\Omega; \delta^{\alpha_+})$; consequently $u \le v$. □

**Proof of Theorem B.** We first notice that since $g \circ (\mathbb{K}_\mu[\nu]) \in L^1_{\text{loc}}(\Omega)$ and $\mathbb{G}_\mu[g \circ (\mathbb{K}_\mu[\nu])] < \infty$, it follows that $g \circ (\mathbb{K}_\mu[\nu]) \in L^1(\Omega; \delta^{\alpha_+})$ due to a similar argument as in the proof of Theorem A. It is easy to see that $\mathbb{K}_\mu[\nu]$ is a subsolution of (3.1). For $\varrho \in (r_1, r_2)$, we look for a supersolution of $v$ of the form

$$v = \varrho \mathbb{K}_\mu[\nu] + c_2 \mathbb{G}_\mu[g \circ (\varrho \mathbb{K}_\mu[\nu])] \tag{3.7}$$

where $c_2$ will be made precise latter on. By (1.15) and (1.17), we obtain

$$v \le \varrho(1 + c_2 c_3 \varrho^{-1} \ell(\varrho)) \mathbb{K}_\mu[\nu].$$

The monotonicity property of $g$ implies

$$g \circ v \le g \circ (\varrho(1 + c_2 c_3 \varrho^{-1} \ell(\varrho)) \mathbb{K}_\mu[\nu]).$$

By (1.15),

$$g \circ v \le \ell(1 + c_2 c_3 \varrho^{-1} \ell(\varrho)) g \circ (\varrho \mathbb{K}_\mu[\nu]). \tag{3.8}$$

In light of (1.16), we deduce

$$g \circ v \le c_2 g \circ (\varrho \mathbb{K}_\mu[\nu]) = -L_\mu v \tag{3.9}$$

This means $v$ is a supersolution of (3.1).

We apply Theorem 3.2 to derive that problem (1.18) admits a minimal solution $\underline{u}_{\varrho\nu}$ satisfying

$$\varrho \mathbb{K}_\mu[\nu] \le \underline{u}_{\varrho\nu} \le \varrho \mathbb{K}_\mu[\nu] + c_2 \mathbb{G}_\mu[g \circ (\mathbb{K}_\mu[\nu])]. \tag{3.10}$$

Estimates (1.19) follows straightforward from (1.17) and (3.10) with $c_4 = 1 + c_2 c_3 \varrho^{-1} \ell(\varrho)$.

We next prove (1.20). Due to (1.13), it is sufficient to prove that for $\nu$-a.e. $z \in \partial\Omega$,

$$\lim_{x \to z} \frac{\mathbb{G}_\mu[g \circ \underline{u}_{\varrho\nu}](x)}{\mathbb{K}_\mu[\nu](x)} = 0 \quad \text{non tangentially.} \tag{3.11}$$

To obtain (3.11), we shall employ Theorem 2.5. Since $\mathbb{K}_\mu[\nu]$ is a positive $L_\mu$-harmonic function satisfying local Harnack inequality, we only need to show that:

(i) $\mathbb{G}_\mu[g \circ \underline{u}_{\varrho\nu}]$ is a positive $L_\mu$-potential.

(ii) $\mathbb{G}_\mu[g \circ \underline{u}_{\varrho\nu}]$ satisfies Harnack inequality.

Since $g \circ \underline{u}_{\varrho\nu} \in L^1(\Omega; \delta^{\alpha_+})$, $\operatorname{tr}^*(\mathbb{G}_\mu[g \circ \underline{u}_{\varrho\nu}])$ and hence (i) follows from Proposition 2.4 (iv). By (1.19), we infer that $\underline{u}_{\varrho\nu}$ satisfies the local Harnack inequality. Since $\underline{u}_{\varrho\nu}$ can be written under the form (1.13), it follows that $\mathbb{G}_\mu[g \circ \underline{u}_{\varrho\nu}]$ satisfies this inequality too. Hence (ii) is verified. By invoking Theorem 2.5, we get (3.11). □

## 4. Power source

In this section, we focus on the equation

$$-L_\mu u = u^q \quad \text{in } \Omega \tag{4.1}$$

4.1. **Subcritical case.** We start with a lemma the proof of which is an adaptation of an idea in [6].

**Lemma 4.1.** *Assume $0 < q < q^*$ and $z \in \partial\Omega$. Then there exists a constant $c_{10} = c_{10}(N, \mu, q, \Omega)$ such that*

$$\mathbb{G}_\mu[K_\mu(\cdot, z)^q](x) \le c_{10}|x - z|^{N+\alpha_+-(N-1-\alpha_-)q} K_\mu(x, z) \quad \forall x \in \Omega. \tag{4.2}$$

**Proof.** By (2.7) and (2.13), there exists a constant $c_{11}$ such that for every $x \in \Omega$,

$$\mathbb{G}_\mu[K_\mu(\cdot, z)^q](x) \le c_{11}\delta(x)^{\alpha_+} \int_\Omega |x - y|^{1+\alpha_--N}|y - z|^{q(1+\alpha_--N)}dy. \tag{4.3}$$

Put

$$\mathcal{D}_1 = \Omega \cap B(x, |x - z|/2), \quad \mathcal{D}_2 = \Omega \cap B(z, |x - z|/2), \quad \mathcal{D}_3 = \Omega \setminus (\mathcal{D}_1 \cup \mathcal{D}_2),$$

and

$$I_i := \int_{\mathcal{D}_i} |x - y|^{1+\alpha_--N}|y - z|^{q(1+\alpha_--N)}dy, \quad i = 1, 2, 3.$$

For every $x \in \mathcal{D}_1$, $|x - z| \le 2|y - z|$, therefore

$$I_1 \le c_{12}|x - z|^{q(1+\alpha_--N)} \int_{\mathcal{D}_1} |x - y|^{1+\alpha_--N}dy \le c'_{12}|x - y|^{1+\alpha_--(N-1-\alpha_-)q}. \tag{4.4}$$

For every $x \in \mathcal{D}_2$, $|x - z| \le 2|x - y|$, hence

$$I_2 \le c_{13}|x - z|^{1+\alpha_--N} \int_{\mathcal{D}_2} |y - z|^{(1+\alpha_--N)q}dy \le c'_{13}|x - y|^{1+\alpha_--(N-1-\alpha_-)q}. \tag{4.5}$$

For every $x \in \mathcal{D}_2$, $|x - y| \ge 3|y - z|$, therefore

$$I_3 \le c_{14} \int_{\mathcal{D}_3} |y - z|^{1+\alpha_--N-(N-1-\alpha_-)q}dy \le c'_{14}|x - y|^{1+\alpha_--(N-1-\alpha_-)q}. \tag{4.6}$$

Combining (4.3)-(4.6), we obtain

$$\mathbb{G}_\mu[K_\mu(\cdot, z)^q](x) \le c_{11}(c'_{12} + c'_{13} + c'_{14})\delta(x)^{\alpha_+}|x - y|^{1+\alpha_--(N-1-\alpha_-)q}. \tag{4.7}$$

Estimate (4.2) follows straightforward from (2.7) and (4.7). □

**Proposition 4.2.** *Assume $0 < q < q^*$ and $\nu$ is a positive finite measure on $\partial\Omega$. Then $\mathbb{K}_\mu[\nu] \in L^q(\Omega; \delta^{\alpha_+})$ and there exists a constant $c_{15}$ depending $N, \mu, q, \Omega$ such that*

$$\mathbb{G}_\mu[\mathbb{K}_\mu[\nu]^q] \le c_{15} \|\nu\|^{q-1}_{\mathfrak{M}(\partial\Omega)} \mathbb{K}_\mu[\nu]. \tag{4.8}$$

**Remark.** We notice that (4.8) was proved in [7] for the case $q > 1$.

**Proof.** We may assume that $\|\nu\|_{\mathfrak{M}(\partial\Omega)} = 1$ (if it is not the case, one can replace $\nu$ by $\nu/\|\nu\|_{\mathfrak{M}(\partial\Omega)}$). We first consider the case $q \ge 1$. From (2.10) and the fact that $L^{q^*}_w(\Omega, \delta^{\alpha_+}) \subset L^q(\Omega, \delta^{\alpha_+})$, we deduce that $\mathbb{K}_\mu[\nu] \in L^q(\Omega; \delta^{\alpha_+})$. It follows from (1.6) and Jensen's inequality that

$$\mathbb{K}_\mu[\nu](x)^q \le \int_{\partial\Omega} K_\mu(x, z)^q d\nu(z).$$

Consequently

$$\mathbb{G}_\mu[\mathbb{K}_\mu[\nu]^q](x) \leq \int_{\partial\Omega} \mathbb{G}_\mu[K_\mu(\cdot,z)^q]d\nu(z).$$

By Lemma 4.1, since $N+\alpha_+-(N-1-\alpha_-)q>0$,

$$\begin{aligned}\mathbb{G}_\mu[\mathbb{K}_\mu[\nu]^q](x) &\leq c_{10}\int_{\partial\Omega}|x-z|^{N+\alpha_+-(N-1-\alpha_-)q}K_\mu(x,z)d\nu(z)\\ &\leq c_{10}(\mathrm{diam}(\Omega))^{N+\alpha_+-(N-1-\alpha_-)q}\mathbb{K}_\mu[\nu](x).\end{aligned}$$

Thus we obtain (4.8).

If $0<q<1$ then

$$\mathbb{G}_\mu[\mathbb{K}_\mu^q[\nu]] \leq \mathbb{G}_\mu[1+\mathbb{K}_\mu[\nu]] = \mathbb{G}_\mu[1]+\mathbb{G}_\mu[\mathbb{K}_\mu[\nu]].$$

From the case $q=1$, we deduce that

$$\mathbb{G}_\mu[\mathbb{K}_\mu^q[\nu]] \leq \mathbb{G}_\mu[1]+c_{15}\mathbb{K}_\mu[\nu]$$

By the estimates $\mathbb{G}_\mu[1]\leq c_{16}\mathbb{K}_\mu[\nu]$, where $c_{16}=c_{16}(N,\mu,\Omega)$, we conclude (4.8). □

**Lemma 4.3.** *Let $f\in L^1(\Omega;\delta^{\alpha_+})$, $f\geq 0$, $\nu\in\mathfrak{M}^+(\partial\Omega)$, $\nu\not\equiv 0$ and $\phi\in C^1([0,\infty))$ be a concave, nondecreasing function such that $\phi(1)\geq 0$ and $\phi'$ is bounded. Let $\varphi$ be a positive function in $L^1_{loc}(\Omega)$ such that $-L_\mu\varphi\geq f$. Then*

$$\phi'\left(\frac{\varphi}{\mathbb{K}_\mu[\nu]}\right)f\in L^1(\Omega;\delta^{\alpha_+}). \tag{4.9}$$

$$-L_\mu\left[\mathbb{K}_\mu[\nu]\phi\left(\frac{\varphi}{\mathbb{K}_\mu[\nu]}\right)\right]\geq\phi'\left(\frac{\varphi}{\mathbb{K}_\mu[\nu]}\right)f \quad \textit{in the weak sense.} \tag{4.10}$$

**Proof.** Since $-L_\mu\varphi\geq f\geq 0$, by Proposition 2.4, there exist $\tau\in\mathfrak{M}^+(\Omega;\delta^{\alpha_+})$ and $\lambda\in\mathfrak{M}^+(\partial\Omega)$ such that

$$\varphi=\mathbb{G}_\mu[f]+\mathbb{G}_\mu[\tau]+\mathbb{K}_\mu[\lambda].$$

Put $\psi=\mathbb{K}_\mu[\nu]$. Let $\{f_n\}$ and $\{\tau_n\}$ be two sequences in $C_c^\infty(\Omega)$ such that $\{f_n\}$ converges to $f$ in $L^1(\Omega;\delta^{\alpha_+})$ and $\{\tau_n\}$ converges to $\tau$ in the weak sense of $\mathfrak{M}^+(\Omega;\delta^{\alpha_+})$. Let $\{\nu_n\}$ and $\{\lambda_n\}$ be two sequences in $C^1(\partial\Omega)$ converging to $\nu$ and $\lambda$ respectively in the weak sense of $\mathfrak{M}^+(\partial\Omega)$. Put $\varphi_n=\mathbb{G}_\mu[f_n]+\mathbb{G}_\mu[\tau_n]+\mathbb{K}_\mu[\lambda_n]$ and $\psi_n=\mathbb{K}_\mu[\nu_n]$. By the bootstrap argument, one can prove that $\varphi_n,\psi_n\in C^2(\Omega)$ for every $n\in\mathbb{N}$. By [16] $\{\mathbb{G}_\mu[f_n]\}$, $\{\mathbb{G}_\mu[\tau_n]\}$, $\{\mathbb{K}_\mu[\lambda_n]\}$ and $\{\mathbb{K}_\mu[\nu_n]\}$ converge to $\mathbb{G}_\mu[f]$, $\mathbb{G}_\mu[\tau]$, $\mathbb{K}_\mu[\lambda]$ and $\mathbb{K}_\mu[\nu]$ respectively in $L^1(\Omega)$. As a consequence, up to a subsequence, $\{\varphi_n\}$ and $\{\psi_n\}$ converge to $\varphi$ and $\psi$ respectively a.e. in $\Omega$. Therefore, for $n$ large enough, $\psi_n>0$.

Due to [9, Lemma 5.3],

$$-\Delta\left[\psi_n\phi\left(\frac{\varphi_n}{\psi_n}\right)\right]\geq\phi'\left(\frac{\varphi_n}{\psi_n}\right)(-\Delta\varphi_n)+\left[\phi\left(\frac{\varphi_n}{\psi_n}\right)-\frac{\varphi_n}{\psi_n}\phi'\left(\frac{\varphi_n}{\psi_n}\right)\right](-\Delta\psi_n).$$

It follows that

$$-L_\mu\left[\psi_n\phi\left(\frac{\varphi_n}{\psi_n}\right)\right]\geq\phi'\left(\frac{\varphi_n}{\psi_n}\right)(-L_\mu\varphi_n)+\left[\phi\left(\frac{\varphi_n}{\psi_n}\right)-\frac{\varphi_n}{\psi_n}\phi'\left(\frac{\varphi_n}{\psi_n}\right)\right](-L_\mu\psi_n).$$

Consequently,

$$-L_\mu\left[\psi_n\phi\left(\frac{\varphi_n}{\psi_n}\right)\right]\geq\phi'\left(\frac{\varphi_n}{\psi_n}\right)f_n. \tag{4.11}$$

Then for every nonnegative function $\zeta \in X(\Omega)$, there holds

$$-\int_\Omega \psi_n \phi\left(\frac{\varphi_n}{\psi_n}\right) L_\mu \zeta dx \geq \int_\Omega \phi'\left(\frac{\varphi_n}{\psi_n}\right) f_n \zeta dx. \tag{4.12}$$

We see that

$$0 \leq \psi_n \phi\left(\frac{\varphi_n}{\psi_n}\right) \leq \psi_n\left(\phi(0) + \phi'(0)\frac{\varphi_n}{\psi_n}\right) = c_{17}(\psi_n + \varphi_n). \tag{4.13}$$

By (2.9) and (2.10), $\{\varphi_n\}$ and $\{\psi_n\}$ are uniformly bounded in $L^p(\Omega;\delta^{-\alpha_-})$ for $p \in (1, \frac{N-\alpha_-}{N-1-\alpha_-})$. Due to Holder inequality, $\{\varphi_n\}$ and $\{\psi_n\}$ are uniformly integrable with respect to $\delta^{\alpha_-}dx$. In view of Vitali theorem $\{\varphi_n\}$ and $\{\psi_n\}$ converge to $\varphi$ and $\psi$ in $L^1(\Omega;\delta^{-\alpha_-})$ respectively. By (4.13) and dominated convergence theorem we deduce that

$$\psi_n \phi\left(\frac{\varphi_n}{\psi_n}\right) \to \psi\phi\left(\frac{\varphi}{\psi}\right) \quad \text{in } L^1(\Omega;\delta^{-\alpha_-}).$$

Due to Fatou lemma, by sending $n \to \infty$ in (4.12), we obtain (4.9) and (4.10). □

**Theorem 4.4.** *Let $q > 1$ and $\nu \in \mathfrak{M}^+(\partial\Omega)$, $\nu \not\equiv 0$. If problem $(D_\nu)$ admits a solution then* (1.22) *holds.*

**Proof.** Since (1.9) admits a solution $u$ then by Theorem A, $u^q \in L^1(\Omega;\delta^{\alpha_+})$ and (1.13) holds. Consequently $\mathbb{K}_\mu[\nu]^q \in L^1(\Omega;\delta^{\alpha_+})$. Now applying Lemma 4.3 with $f = u^q$, $\varphi = u$ and

$$\phi(s) = \begin{cases} \dfrac{1-s^{1-q}}{q-1}, & \text{if } s \geq 1, \\ s-1 & \text{if } 0 \leq s < 1 \end{cases}$$

we obtain the following estimate in the weak sense

$$-L_\mu\left[\mathbb{K}_\mu[\nu]\phi\left(\frac{u}{\mathbb{K}_\mu[\nu]}\right)\right] \geq \left(\frac{u}{\mathbb{K}_\mu[\nu]}\right)^{-q} u^q \geq \mathbb{K}_\mu[\nu]^q. \tag{4.14}$$

Put

$$\Psi := \mathbb{K}_\mu[\nu]\phi\left(\frac{u}{\mathbb{K}_\mu[\nu]}\right) \quad \text{and} \quad \tilde{\Psi} := \mathbb{G}_\mu[\mathbb{K}_\mu[\nu]^q].$$

Then $\Psi$ is an $L_\mu$-superharmonic function and by Proposition 2.4, $\Psi$ admits a nonnegative normalized boundary trace. By Kato lemma (see [20]), $(\tilde{\Psi} - \Psi)_+$ is an $L_\mu$-subharmonic function and $\operatorname{tr}^*((\tilde{\Psi} - \Psi)_+) = 0$. It follows that $(\tilde{\Psi} - \Psi)_+ = 0$ and hence $\tilde{\Psi} \leq \Psi$ in $\Omega$. This means

$$\mathbb{G}_\mu[\mathbb{K}_\mu[\nu]^q] \leq \mathbb{K}_\mu[\nu]\phi\left(\frac{u}{\mathbb{K}_\mu[\nu]}\right) \leq \frac{1}{q-1}\mathbb{K}_\mu[\nu].$$

□

**Proof of Proposition C.** This theorem follows straightforward from Lemma 4.1 and Theorem 4.4. □

**Proposition 4.5.** *Assume $0 < q < q^*$, $q \neq 1$ and $\nu \in \mathfrak{M}^+(\partial\Omega)$ such that $\|\nu\|_{\mathfrak{M}(\partial\Omega)} = 1$.*

(i) *If $q > 1$ then there exists a positive number $\varrho_0 > 0$ depending on $N, \mu, q, \Omega$ such that for every $\varrho \in (0, \varrho_0)$ problem $(D_{\varrho\nu})$ admits a minimal solution $\underline{u}_{\varrho\nu}$.*

(ii) *If $q \in (0,1)$ then for every $\varrho > 0$ problem $(D_{\varrho\nu})$ admits a minimal solution $\underline{u}_{\varrho\nu}$.*

*In any case $\underline{u}_{\varrho\nu}$ satisfies* (1.19) *and* (1.20).

**Proof.** We shall apply Theorem B to deduce the existence of a solution of $(D_{\varrho\nu})$. One can verify that the functions $g(x,s) = s^q$ and $\ell(s) = s^q$ with $q > 0$ satisfy (1.15). From Proposition 4.2 we deduce that condition (1.17) is fulfilled with the constant $c_{15}$. For such $g$ and $\ell$, condition (1.16) is valid if one can find a positive constant $c_{18}$ such that

$$1 + c_{18}c_{15}\varrho^{q-1} \le c_{18}^{\frac{1}{q}}. \tag{4.15}$$

If $q > 1$ then there exists $\varrho_0 = \varrho_0(q, c_{15})$ and $c_{18} = c_{18}(q)$ such that (4.15) holds true. If $q < 1$ then for every $\varrho \in [1,\infty)$ one can choose $c_{18} = c_{18}(q, c_{15})$ such that (4.15) holds. Hence, by Theorem B, there exists a minimal solution $\underline{u}_{\varrho\nu}$ of problem $(D_{\varrho\nu})$ if $q > 1$ and $\varrho \in (0, \varrho_0)$ or if $q < 1$ and $\varrho \in [1,\infty)$.

We next consider the case $q < 1$ and $0 < \varrho < 1$. Let $v$ be a solution of the problem $(D_\nu)$ then $v \ge \varrho\mathbb{K}_\mu[\nu]$. Due to Theorem 3.2, we deduce the existence of a minimal solution $\underline{u}_{\varrho\nu}$ of problem $(D_{\varrho\nu})$ satisfying $\varrho\mathbb{K}_\mu[\nu] \le \underline{u}_{\varrho\nu} \le v$.

Thus, if $q < 1$, for any $\varrho > 0$, $\underline{u}_{\varrho\nu}$ satisfies (1.19) and (1.20). □

**Lemma 4.6.** *Let $0 < q \ne 1$ and $\nu \in \mathfrak{M}^+(\partial\Omega)$. If $u$ is a solution of $(D_\nu)$ then there is a constant $c_{19} = c_{19}(N,\mu,q,\Omega)$ such that*

$$\|u\|_{L^1(\Omega;\delta^{-\alpha_-})} + \|u\|_{L^q(\Omega;\delta^{\alpha_+})} \le c_{19}(1 + \|\nu\|_{\mathfrak{M}(\partial\Omega)}). \tag{4.16}$$

**Proof.** Indeed, by taking $\zeta = \varphi_{\mu,1}$ in the formulation satisfied by $u$, we obtain

$$\lambda_{\mu,1}\int_\Omega u\varphi_{\mu,1}dx = \int_\Omega u^q\varphi_{\mu,1}dx + \lambda_{\mu,1}\int_\Omega \mathbb{K}_\mu[\nu]\varphi_{\mu,1}dx \quad \forall\zeta \in X(\Omega). \tag{4.17}$$

**Case 1:** $q > 1$**.** By Young inequality, we get

$$\int_\Omega u\varphi_{\mu,1}dx \le (2\lambda_{\mu,1})^{-1}\int_\Omega u^q\varphi_{\mu,1}dx + (2\lambda_{\mu,1})^{\frac{1}{q-1}}\int_\Omega \varphi_{\mu,1}dx. \tag{4.18}$$

By (4.17) and (4.18), we obtain

$$\int_\Omega u^q\varphi_{\mu,1}dx + 2\lambda_{\mu,1}\int_\Omega \mathbb{K}_\mu[\nu]\varphi_{\mu,1}dx \le (2\lambda_{\mu,1})^{\frac{q}{q-1}}\int_\Omega \varphi_{\mu,1}dx. \tag{4.19}$$

Since the second term on the left hand-side of (4.19) is nonnegative and by (1.8), we get

$$\|u\|^q_{L^q(\Omega;\delta^{\alpha_+})} \le c_1^2(2\lambda_{\mu,1})^{\frac{q}{q-1}}\int_\Omega \delta^{\alpha_+}dx \le c_{20}. \tag{4.20}$$

On the other hand, we derive from (1.13), (2.9) and (2.10) that

$$\|u\|_{L^1(\Omega;\delta^{-\alpha_-})} \le c_{21}(\|u\|_{L^q(\Omega;\delta^{\alpha_+})} + \|\nu\|_{\mathfrak{M}(\partial\Omega)}). \tag{4.21}$$

Combining (4.20) and (4.21), we obtain (4.16).

**Case 2:** $q \in (0,1)$**.** By Young inequality, we have

$$\int_\Omega u^q\varphi_{\mu,1}dx \le \frac{\lambda}{2}\int_\Omega u\varphi_{\mu,1}dx + (2\lambda_{\mu,1}^{-1})^{\frac{q}{1-q}}\int_\Omega \varphi_{\mu,1}dx. \tag{4.22}$$

Consequently,

$$\int_\Omega u\varphi_{\mu,1}dx \le (2\lambda_{\mu,1}^{-1})^{\frac{1}{1-q}}\int_\Omega \varphi_{\mu,1}dx + 2\int_\Omega \mathbb{K}_\mu[\nu]\varphi_{\mu,1}dx.$$

Therefore

$$\|u\|_{L^1(\Omega;\delta^{-\alpha_-})} \le c_{22}(1 + \|\nu\|_{\mathfrak{M}(\partial\Omega)}). \tag{4.23}$$

Combining (4.25) and (4.23) leads to (4.16). □

**Theorem 4.7.** *Assume $q \in (1, q^*)$. Let $\nu \in \mathfrak{M}^+(\partial\Omega)$ with $\|\nu\|_{\mathfrak{M}(\partial\Omega)} = 1$. Then there exists a threshold value $\varrho^* \in \mathbb{R}_+$ for $(D_{\varrho\nu})$ such that*

(i) *If $\varrho \in (0, \varrho^*]$ then problem $(D_{\varrho\nu})$ admits a minimal solution $\underline{u}_{\varrho\nu}$. If $\varrho \in (0, \varrho^*)$ then $\underline{u}_{\varrho\nu}$ satisfies* (1.19) *and* (1.20)*.*

*Moreover $\{\underline{u}_{\varrho\nu}\}$ is an increasing sequence which converges, as $\varrho \to \varrho^*$, to the minimal solution $\underline{u}_{\varrho^*\nu}$ of $(D_{\varrho^*\nu})$ in $L^1(\Omega; \delta^{-\alpha_-})$ and in $L^1(\Omega; \delta^{\alpha_+})$.*

(ii) *If $\varrho > \varrho^*$ then there exists no solution of $(D_{\varrho\nu})$.*

**Proof.** Put

$$\mathcal{A} := \{\varrho > 0 : (D_{\varrho\nu}) \text{ admits a solution}\} \quad \text{and} \quad \varrho^* = \sup \mathcal{A}.$$

By Proposition 4.5, $(D_{\varrho\nu})$ admits a solution for $\varrho > 0$ small, therefore $\mathcal{A} \neq \emptyset$. Moreover, from Theorem 4.4, we deduce that $\varrho^*$ is finite.

We shall show that $(0, \varrho^*) \subset \mathcal{A}$. For this purpose, we have to show that if $0 < \varrho < \varrho'$ and $\mathcal{A} \ni \varrho' < \varrho^*$ then $\varrho \in \mathcal{A}$. Since $\varrho' \in \mathcal{A}$, due to Theorem 4.4, there exists a minimal solution $\underline{u}_{\varrho'\nu}$ of $(D_{\varrho'\nu})$ which is greater than $\varrho\mathbb{K}_\mu[\nu]$. By Theorem 3.2, problem $(D_{\varrho\nu})$ admits a minimal solution $\underline{u}_{\varrho\nu}$, i.e. $\varrho \in \mathcal{A}$.

Next we prove that $\varrho^* \in \mathcal{A}$, namely problem $(D_{\varrho^*\nu})$ admits a solution. Let $\{\varrho_n\}$ be an increasing sequence converging to $\varrho^*$. For each $n$, let $u_{\varrho_n\nu}$ be a solution of $(D_{\varrho_n\nu})$. Then $u_{\varrho_n\nu} \in L^1(\Omega; \delta^{-\alpha_-}) \cap L^q(\Omega; \delta^{\alpha_+})$ and it satisfies the formula

$$-\int_\Omega u_{\varrho_n\nu} L_\mu \zeta dx = \int_\Omega u_{\varrho_n\nu}^q \zeta dx - \varrho_n \int_\Omega \mathbb{K}_\mu[\nu] L_\mu \zeta dx \quad \forall \zeta \in X(\Omega). \tag{4.24}$$

It follows from Lemma 4.6 that the sequence $\{u_{\varrho_n\nu}^q\}$ is uniformly bounded in $L^1(\Omega; \delta^{\alpha_+})$ and hence by local regularity for elliptic equations [18] there exists a subsequence, still denoted by the same notation, such that $\{u_{\varrho_n\nu}\}$ converges a.e. to a function $u_{\varrho^*\nu}$. From Theorem A, there holds

$$u_{\varrho_n\nu} = \mathbb{G}_\mu[u_{\varrho_n\nu}^q] + \varrho_n \mathbb{K}_\mu[\nu]. \tag{4.25}$$

Thanks to Proposition 2.1, $\{u_{\varrho_n\nu}\}$ is uniformly bounded in $L^{q_1}(\Omega; \delta^{-\alpha_-})$ and in $L^{q_2}(\Omega; \delta^{\alpha_+})$ where $1 < q_1 < \frac{N-\alpha_-}{N-1-\alpha_-}$ and $q < q_2 < q^*$. We invoke Holder inequality to infer that $\{u_{\varrho_n\nu}\}$ and $\{u_{\varrho_n\nu}^q\}$ are uniformly integrable with respect to $\delta^{-\alpha_-}dx$ and $\delta^{\alpha_+}dx$ respectively. As a consequence, $\{u_{\varrho_n\nu}\}$ converges to $u_{\varrho^*\nu}$ in $L^1(\Omega; \delta^{-\alpha_-})$ and $\{u_{\varrho_n\nu}^q\}$ converges to $u_{\varrho^*\nu}^q$ in $L^1(\Omega; \delta^{\alpha_+})$. Letting $n \to \infty$ in (4.24) implies

$$-\int_\Omega u_{\varrho^*\nu} L_\mu \zeta dx = \int_\Omega u_{\varrho^*\nu}^q \zeta dx - \varrho^* \int_\Omega \mathbb{K}_\mu[\nu] L_\mu \zeta dx \quad \forall \zeta \in X(\Omega). \tag{4.26}$$

We infer from Theorem A that $u_{\varrho^*\nu}$ is a solution of $(D_{\varrho^*\nu})$.

Notice that, in light of Theorem 3.2 and the above argument, one can prove that $\{\underline{u}_{\varrho\nu}\}$ is an increasing sequence converging to the minimal solution $\underline{u}_{\varrho^*\nu}$ of $(D_{\varrho^*\nu})$ in $L^1(\Omega; \delta^{-\alpha_-})$ and in $L^1(\Omega; \delta^{\alpha_+})$.

We next show that for each $\varrho \in (0, \varrho^*)$, there exists a minimal solution $\underline{u}_{\varrho\nu}$ of $(D_{\varrho\nu})$ which satisfies (1.19). Take $\varrho' = \frac{\varrho+\varrho^*}{2}$ and let $u_{\varrho'\nu}$ be a solution of $(D_{\varrho'\nu})$. We next apply (4.10) with $\nu$ replaced by $\varrho'\nu$, $\varphi = u_{\varrho'\nu}$, $f = u_{\varrho'\nu}^q$ and

$$\phi(s) = \begin{cases} s(1 + \varepsilon s^{q-1})^{-\frac{1}{q-1}}, & \text{if } s \geq 1, \\ \left(\dfrac{\varrho}{\varrho'}\right)^q s + \left(\dfrac{\varrho}{\varrho'}\right) - \left(\dfrac{\varrho}{\varrho'}\right)^q & \text{if } 0 \leq s \leq 1 \end{cases} \quad \text{with} \quad \varepsilon = \left(\frac{\varrho'}{\varrho}\right)^{q-1} - 1.$$

We get

$$-L_\mu\left(\mathbb{K}_\mu[\varrho'\nu]\phi\left(\frac{u_{\varrho'\nu}}{\mathbb{K}_\mu[\varrho'\nu]}\right)\right)\geq \phi'\left(\frac{u_{\varrho'\nu}}{\mathbb{K}_\mu[\varrho'\nu]}\right)u_{\varrho'\nu}^q=\left(\mathbb{K}_\mu[\varrho'\nu]\phi\left(\frac{u_{\varrho'\nu}}{\mathbb{K}_\mu[\varrho'\nu]}\right)\right)^q.$$

Therefore

$$\Psi=\mathbb{K}_\mu[\varrho'\nu]\phi\left(\frac{u_{\varrho'\nu}}{\mathbb{K}_\mu[\varrho'\nu]}\right)$$

is a supersolution of (4.1). Moreover $\Psi\geq\varrho\mathbb{K}_\mu[\nu]$. By Theorem 3.2 there exists a minimal solution $\underline{u}_{\varrho\nu}$ of $(D_{\varrho\nu})$ such that

$$\varrho\mathbb{K}_\mu[\nu]\leq\underline{u}_{\varrho\nu}\leq\Psi.$$

This implies

$$\varrho\mathbb{K}_\mu[\nu]\leq\underline{u}_{\varrho\nu}\leq\varepsilon^{-\frac{1}{q-1}}\varrho'\mathbb{K}_\mu[\nu].$$

Therefore we get (1.19) with $c_4=\varrho^{-1}\varrho'\varepsilon^{-\frac{1}{q-1}}$. Finally, (1.20) can be obtained by a similar argument as in the proof of Theorem B. □

**Proof of Theorem D.** Part (1) follows from Theorem 4.7. Part (2) follows from Proposition 4.5 (ii).

If $\nu=\varrho\delta_z$, by (1.13) and (1.19), we obtain

$$\varrho\leq\frac{\underline{u}_{\varrho\delta_z}(x)}{K_\mu(x,z)}\leq\varrho+c_{23}\frac{\mathbb{G}_\mu[K_\mu(\cdot,z)^q](x)}{K_\mu(x,z)}\tag{4.27}$$

Since $q<q^*$, it follows from Lemma 4.1 that

$$\lim_{x\to z}\frac{\mathbb{G}_\mu[K_\mu(\cdot,z)^q](x)}{K_\mu(x,z)}=0.$$

Thus, by (4.27), we conclude (1.23). □

**Proof of Proposition E.** If $q>1$, assumption (1.24) guarantees the existence of a solution $u_{\nu_n}$ of $(D_{\nu_n})$. Moreover, since $\{\nu_n\}$ converges weakly to $\nu$, it follows that $\|\nu\|_{\mathfrak{M}(\partial\Omega)}\leq\varrho^*$. Due to Lemma 4.6, the sequence $\{u_{\nu_n}\}$ is uniformly bounded in $L^q(\Omega;\delta^{\alpha_+})$. Employing a similar argument as in the proof of Theorem 4.7, we obtain the convergence in $L^1(\Omega;\delta^{-\alpha_-})$ and in $L^q(\Omega;\delta^{\alpha_+})$.

If $q\in(0,1)$, due to Lemma 4.6, we obtain the convergence in $L^1(\Omega;\delta^{-\alpha_-})$. □

We next consider the case $q=1$.

**Lemma 4.8.** *Let $\kappa>0$ and $u$ be a positive solution of*

$$-L_\mu u=\kappa u\quad in\ \Omega.\tag{4.28}$$

*Then $u$ satisfies the Harnack inequality; i.e. for every $a\in(0,1)$, there exists a constant $c_{24}=c_{24}(N,\mu,q,\Omega)$ such that for every $x\in\Omega$,*

$$\sup_{B(x,a\delta(x))}u\leq c_{24}\inf_{B(x,a\delta(x))}u.\tag{4.29}$$

**Proof.** Equation (4.28) can be written as follows

$$-\Delta u=\left(\frac{\mu}{\delta^2}+\kappa\right)u\qquad\text{in }\Omega.\tag{4.30}$$

Take arbitrarily $a\in(0,1)$ and $x_0\in\Omega$. Put $d=\frac{a+1}{2}\delta(x_0)$ and $M:=\max_{B(x_0,d)}u$. Put

$$y_0=d^{-1}x_0,\quad\Omega_d=d^{-1}\Omega,\quad\delta_d(y)=\operatorname{dist}(y,\partial\Omega_d)\text{ with }y\in\Omega_d.$$

We define

$$v_d(y) := M^{-1}u(dy), \quad \forall y \in \Omega_d.$$

Clearly, $\max_{B(y_0;1)} v_d = 1$ and due to (4.30) we deduce that $v_d$ is a solution of

$$-\Delta v_d = V v_d \qquad \text{in } \Omega_d. \tag{4.31}$$

where

$$V(y) := \frac{\mu}{\delta_d(y)^2} + d^2\kappa.$$

One can find a positive constant $c_{25}$ such that $V(y) \le c_{25}\delta_d(y)^{-2}$ for every $y \in \Omega_d$. Notice that $B(y_0;1) \subset \Omega_d$ and for every $y \in B(y_0;1)$, there holds

$$\delta_d(y) \ge \frac{1-a}{1+a}.$$

Hence $0 \le V \le c_{26}$ in $B(y_0;1)$ where $c_{26} = c_{26}(a,\mu)$ . By applying Harnack inequality, we deduce that there is a constant $c_{27} = c_{27}(a,\mu,N,\Omega)$ such that

$$\sup_{B(y_0;\frac{2a}{a+1})} v_d \le c_{27} \inf_{B(y_0;\frac{2a}{a+1})} v_d.$$

Thus we obtain (4.29). □

**Proof of Theorem F.** Put

$$\kappa_0 = \min\{1, \|\mathbb{G}_\mu[1]\|^{-1}_{L^\infty(\Omega)}\}.$$

**Claim 1.** For any $\kappa \in (0,\kappa_0)$ there exists a minimal solution $\underline{u}_{\kappa,\nu}$ of problem $(E^\kappa_\nu)$.

Fix $q \in (1,q^*)$ such that $\kappa < (q\,\|\mathbb{G}_\mu[1]\|_{L^\infty(\Omega)})^{-1}$ and let $\varrho^*$ be the threshold value for $(D_{\varrho\nu})$. Put $\varrho = \|\nu\|_{\mathfrak{M}(\partial\Omega)} > 0$ and $\tilde\nu = \nu/\varrho$.

We first assume that $\varrho \in (0,\varrho^*)$ and let $\underline{u}_\nu$ be the minimal solution of the problem $(D_\nu)$. By Young inequality, we get

$$u_\nu^q + 1 \ge u_\nu + \frac{1}{q} \ge \kappa(u_\nu + \mathbb{G}_\mu[1]).$$

It follows that

$$-L_\mu(u_\nu + \mathbb{G}_\mu[1]) \ge \kappa(u_\nu + \mathbb{G}_\mu[1]).$$

Therefore $u_\nu + \mathbb{G}_\mu[1]$ is a super solution of the equation

$$-L_\mu u = \kappa u. \tag{4.32}$$

Clearly $\mathbb{K}_\mu[\nu]$ is a subsolution of (4.32). By Theorem 3.2 there is a minimal solution $\underline{u}_{\kappa,\nu}$ of $(E^\kappa_\nu)$ which satisfies $\mathbb{K}_\mu[\nu] \le \underline{u}_{\kappa,\nu} \le \underline{u}_\nu + \mathbb{G}_\mu[1]$. Since $\mathbb{G}_\mu[1] \le c_{16}\mathbb{K}_\mu[\nu]$, we infer that $\underline{u}_{\kappa,\nu}$ satisfies (1.19) and (1.20).

If $\varrho \ge \varrho^*$ then there exists $m > 0$ such that $\varrho/m \in (0,\varrho^*)$. Let $\underline{u}_{\kappa,\frac{\nu}{m}}$ be a solution of $(E^\kappa_{\frac{\nu}{m}})$. Put $\underline{u}_{\kappa,\nu} = m\underline{u}_{\kappa,\frac{\nu}{m}}$ then by the linearity, we deduce that $\underline{u}_{\kappa,\nu}$ is the minimal solution of $(E^\kappa_\nu)$ with satisfies (1.19) and (1.20).

**Claim 2.** There exists a number $\kappa^* \in (0,\lambda_{\mu,1}]$ such that

(i) If $\kappa \in (0,\kappa^*)$ then $(E^\kappa_\nu)$ admits a solution.

(ii) If $\kappa > \kappa^*$ then $(E^\kappa_\nu)$ admits no solution.

Put $\mathcal{B} := \{\kappa > 0 : (E^\kappa_\nu) \text{ admits a solution }\}$ and denote $\kappa^* := \sup\mathcal{B}$. We shall show that $(0,\kappa^*) \subset \mathcal{B}$. Take $\kappa' \in \mathcal{B}$ and let $u_{\kappa',\nu}$ be the minimal solution of $(E^{\kappa'}_\nu)$. For any $\kappa \in (0,\kappa')$, $u_{\kappa',\nu}$ and $\mathbb{K}_\mu[\nu]$ are respectively super and sub solutions of $(E^\kappa_\nu)$ such that $\mathbb{K}_\mu[\nu] \le u_{\kappa',\nu}$. Then by Theorem 3.2 there exists a minimal solution $\underline{u}_{\kappa,\nu}$ of $(E_\nu^\kappa)$ satisfying $\mathbb{K}_\mu[\nu] \leq \underline{u}_{\kappa,\nu} \leq u_{\kappa',\nu}$ in $\Omega$. Hence $\kappa \in \mathcal{B}$.

By Lemma 4.8, $\mathbb{G}_\mu[\underline{u}_{\kappa,\nu}]$ satisfies local Harnack inequality. Hence, we deduce from Theorem 2.5 that, for $\nu$-a.e. $z \in \partial\omega$, there holds

$$\lim_{x \to z} \frac{\mathbb{G}_\mu[\underline{u}_{\kappa,\nu}](x)}{\mathbb{K}_\mu[\nu](x)} = 0.$$

Consequently, (1.20) remains valid with $\underline{u}_{\varrho\nu}$ replaced by $\underline{u}_{\kappa,\nu}$.

Now let $\nu \in \mathfrak{M}^+(\partial\Omega)$, $\kappa \in \mathcal{B}$ and denote by $u_{\kappa,\nu}$ a solution of $(E_\nu^\kappa)$. Then by Theorem A,

$$-\int_\Omega u_{\kappa,\nu} L_\mu \zeta dx = \kappa \int_\Omega u_{\kappa,\nu} \zeta dx - \int_\Omega \mathbb{K}_\mu[\nu] L_\mu \zeta dx \quad \forall \zeta \in X(\Omega).$$

Taking $\zeta = \varphi_{\mu,1}$, we obtain

$$\lambda_{\mu,1} \int_\Omega u_{\kappa,\nu} \varphi_{\mu,1} dx = \kappa \int_\Omega u_{\kappa,\nu} \varphi_{\mu,1} dx + \lambda_{\mu,1} \int_\Omega \mathbb{K}_\mu[\nu] \varphi_{\mu,1} dx, \tag{4.33}$$

which implies that $\kappa < \lambda_{\mu,1}$. Consequently, $\kappa^* \leq \lambda_{\mu,1}$.

We show that $\lambda_{\mu,1} \notin \mathcal{B}$ by contradiction. Indeed, suppose that there exists $\nu \in \mathfrak{M}^+(\partial\Omega)$ such that the problem

$$-L_\mu u = \lambda_{\mu,1} u \quad \text{in } \Omega, \qquad \operatorname{tr}^*(u) = \nu \tag{4.34}$$

admits a solution $\hat{u}$. Take $\varphi_{\mu,1}$ as a test function in the weak formulation satisfied by $\hat{u}$, we deduce $\nu \equiv 0$, which is a contradiction.

Now let $\kappa \in (0, \kappa^*)$ and assume $\{\nu_n\}$ is a sequence of measures in $\mathfrak{M}^+(\partial\Omega)$ which converges weakly to $\nu \in \mathfrak{M}^+(\partial\Omega)$. Let $u_{\kappa,\nu_n}$ be a solution of $(E_{\nu_n}^\kappa)$. By (4.33), we deduce

$$\|u_{\kappa,\nu_n}\|_{L^1(\Omega;\delta^{\alpha_+})} \leq c_{28} (\lambda_{\mu,1} - \kappa)^{-1} \|\nu_n\|_{\mathfrak{M}(\partial\Omega)} \leq c_{29} (\lambda_{\mu,1} - \kappa)^{-1} \|\nu\|_{\mathfrak{M}(\partial\Omega)}.$$

By a similar argument as in the proof of Theorem 4.7, we deduce that, up to a subsequence, $\{u_{\nu,\nu_n}\}$ converges to a solution $u_{\kappa,\nu}$ of $(E_\nu^\kappa)$ in $L^1(\Omega, \delta^{-\alpha_-})$. □

**Remark.** (i) If $\kappa > 0$ small then $\underline{u}_{\kappa,\nu}$ satisfies (1.19). Moreover, if $\nu = \varrho\delta_z$ with $\varrho > 0$, $z \in \partial\Omega$ then $\underline{u}_{\kappa,\varrho\delta_z}$ satisfies (1.23).

(ii) A question remains open: "Is $\kappa^* = \lambda_{\mu,1}$?" In case that this equality holds true then $(E_\nu^{\kappa^*})$ admits no solution. Otherwise, if $\kappa^* < \lambda_{\mu,1}$, $(E_\nu^{\kappa^*})$ admits a solution.

## 5. Subcriticality and sublinearity

In this section, we assume that $(g \circ u)(x) = \delta(x)^\gamma \tilde{g}(u(x))$ where $\gamma > -1 - \alpha_+$ and $\tilde{g} \in C(\mathbb{R}_+)$, $\tilde{g}(0) = 0$. The proof of Theorems H and I is an adaptation of the idea in [11]. A distinct feature of this approach is that convexity and monotonicity hypotheses of $g$ can be relaxed while these properties are crucial in other methods.

5.1. **Subcriticality.** Let $\{g_n\}$ be a sequence of $C^1$ nonnegative functions defined on $\mathbb{R}_+$ such that

$$g_n(0) = \tilde{g}(0),\ g_n \leq g_{n+1} \leq \tilde{g},\ \sup_{\mathbb{R}_+} g_n = n \text{ and } \lim_{n \to \infty} \|g_n - \tilde{g}\|_{L^\infty_{loc}(\mathbb{R}_+)} = 0. \tag{5.1}$$

Put

$$\tilde{\gamma} = \min\{\alpha_+ + \gamma, -\alpha_-\} > -1. \tag{5.2}$$

In preparation for proving Theorem G, we establish the following lemma:

**Lemma 5.1.** *Let $\nu \in \mathfrak{M}^+(\partial\Omega)$ such that $\|\nu\|_{\mathfrak{M}(\partial\Omega)} = 1$ and $\{g_n\} \subset C^1(\mathbb{R}_+)$ be a sequence satisfying (5.1). Assume (1.27) and (1.28) are satisfied. Then there exist $\bar\lambda$, $\theta_0 > 0$ and $\varrho_0 > 0$ depending on $\Lambda_0$, $\Lambda_1$, $N$, $\mu$, $\gamma$ and $q_1$ such that for every $\theta \in (0,\theta_0)$ and $\varrho \in (0,\varrho_0)$ the following problem*

$$-L_\mu v = \delta^\gamma g_n(v + \varrho\mathbb{K}_\mu[\nu]) \quad \text{in } \Omega, \qquad \operatorname{tr}^*(v) = 0 \tag{5.3}$$

*admits a nonnegative $v_n \in L_w^{q_\gamma^*}(\Omega;\delta^{\alpha_+ + \gamma}) \cap L^{q_1}(\Omega;\delta^{-\alpha_-})$ satisfying*

$$\|v_n\|_{L_w^{q_\gamma^*}(\Omega;\delta^{\alpha_++\gamma})} + \|v_n\|_{L^{q_1}(\Omega;\delta^{\tilde\gamma})} \le \bar\lambda. \tag{5.4}$$

**Proof.** We shall use Schauder fixed point theorem to show the existence of positive solutions of (5.3). For $n \in \mathbb{N}$, define the operator $\mathbb{S}_n$ by

$$\mathbb{S}_n(v) = \mathbb{G}_\mu[\delta^\gamma g_n(v + \varrho\mathbb{K}_\mu[\nu])], \quad \forall v \in L^1_+(\Omega). \tag{5.5}$$

Set

$$\begin{aligned} M_1(v) &= \|v\|_{L_w^{q_\gamma^*}(\Omega;\delta^{\alpha_++\gamma})}, \quad \forall v \in L_w^{q_\gamma^*}(\Omega;\delta^{\alpha_++\gamma}) \\ M_2(v) &= \|v\|_{L^{q_1}(\Omega;\delta^{\tilde\gamma})}, \quad \forall v \in L^{q_1}(\Omega;\delta^{\tilde\gamma}) \\ M(v) &= M_1(v) + M_2(v), \quad \forall v \in L_w^{q_\gamma^*}(\Omega;\delta^{\alpha_++\gamma}) \cap L^{q_1}(\Omega;\delta^{\tilde\gamma}). \end{aligned} \tag{5.6}$$

**Step 1**: Estimate on $g_n(v + \varrho\mathbb{K}_\mu[\nu])$ in $L^1(\Omega;\delta^{\alpha_++\gamma})$ for $v \in L_w^{q_\gamma^*}(\Omega;\delta^{\alpha_++\gamma}) \cap L^{q_1}(\Omega;\delta^{\tilde\gamma})$.

For $\lambda > 0$, set $A_\lambda = \{x \in \Omega : v + \varrho\mathbb{K}_\mu[\nu] > \lambda\}$ and $a(\lambda) = \int_{A_\lambda} \delta^{\alpha_++\gamma}dx$. We write

$$\begin{aligned} \|g_n(v+\varrho\mathbb{K}_\mu[\nu])\|_{L^1(\Omega;\delta^{\alpha_++\gamma})} &= \int_{A_1} g_n(v+\varrho\mathbb{K}_\mu[\nu])\delta^{\alpha_++\gamma}dx + \int_{A_1^c} g_n(v+\varrho\mathbb{K}_\mu[\nu])\delta^{\alpha_++\gamma}dx \\ &=: I + II. \end{aligned} \tag{5.7}$$

We first estimate $I$ from above. We see that

$$I = a(1)g_n(1) + \int_1^\infty a(s)dg_n(s).$$

Since (1.27) holds, it was proved in [11, Lemma 3.1] that there exists an increasing sequence of real positive number $\{\ell_j\}$ such that

$$\lim_{j\to\infty} \ell_j = \infty \quad \text{and} \quad \lim_{j\to\infty} \ell_j^{-q_\gamma^*}\tilde g(\ell_j) = 0. \tag{5.8}$$

Consequently,

$$\lim_{j\to\infty} \ell_j^{-q_\gamma^*} g_n(\ell_j) = 0, \quad \forall n \in \mathbb{N}. \tag{5.9}$$

Observe that

$$\int_1^\infty a(s)dg_n(s) = \lim_{j\to\infty}\int_\lambda^{\ell_j} a(s)dg_n(s).$$

On the other hand, by (2.5) one gets, for every $s > 0$,

$$a(s) \le \|v + \varrho\mathbb{K}_\mu[\nu]\|^{q_\gamma^*}_{L_w^{q_\gamma^*}(\Omega;\delta^{\alpha_++\gamma})} s^{-q_\gamma^*} \le c_{31}(M_1(v) + c_{32}\varrho)^{q_\gamma^*}s^{-q_\gamma^*} \tag{5.10}$$

where $c_i = c_i(N,\mu,\Omega)$ with $i = 31, 32$. Using (5.10), we obtain

$$\begin{aligned}
a(1)g_n(1) &+ \int_1^{\ell_j} a(s)dg_n(s) \\
&\leq c_{31}(M_1(v)+c_{32}\varrho)^{q_\gamma^*} g_n(1) + c_{31}(M_1(v)+c_{32}\varrho)^{q_\gamma^*}\int_1^{\ell_j} s^{-q_\gamma^*}dg_n(s) \\
&\leq c_{31}(M_1(v)+c_{32}\varrho)^{q_\gamma^*}\ell_j^{-q_\gamma^*} g_n(\ell_j) + c_{31}q_\gamma^*(M_1(v)+c_{32}\varrho)^{q_\gamma^*}\int_1^{\ell_j} s^{-1-q_\gamma^*}g_n(s)ds.
\end{aligned}$$

By virtue of (5.8), letting $j \to \infty$ yields

$$I \leq c_{31}q_\gamma^*(M_1(v)+c_{32}\varrho)^{q_\gamma^*}\int_1^{\infty} s^{-1-q_\gamma^*}g_n(s)ds \leq c_{33}\Lambda_0 M_1(v)^{q_\gamma^*} + c_{33}\Lambda_0\varrho^{q_\gamma^*} \tag{5.11}$$

where $c_{33} = c_{33}(N,\mu,\Omega)$.

To handle the remaining term $II$, without lost of generality, we assume $q_1 \in (1, \frac{N+\tilde{\gamma}}{N-1-\alpha_-})$. Since $\tilde{g}$ satisfies condition (1.28) and $g_n \leq \tilde{g}$, it follows that $g_n$ satisfies this condition too. Hence

$$\begin{aligned}
II &\leq \Lambda_1\int_{A_1^c}(v+\varrho\mathbb{K}_\mu[\nu])^{q_1}\delta^{\alpha_++\gamma}dx + \theta\int_{A_1^c}\delta^{\alpha_++\gamma}dx \\
&\leq \Lambda_1 c_{34}\int_\Omega v^{q_1}\delta^{\alpha_++\gamma}dx + \Lambda_1 c_{34}\varrho^{q_1} + c_{34}\theta \\
&\leq \Lambda_1 c_{35}M_2(v)^{q_1} + \Lambda_1 c_{34}\varrho^{q_1} + c_{34}\theta
\end{aligned} \tag{5.12}$$

where $c_i = c_i(N,\mu,\Omega)$, $i = 33, 34$.

Combining (5.7), (5.11) and (5.12) yields

$$\|g_n(v+\varrho\mathbb{K}_\mu[\nu])\|_{L^1(\Omega;\delta^{\alpha_++\gamma})} \leq c_{33}\Lambda_0 M_1(v)^{q_\gamma^*} + c_{35}\Lambda_1 M_2(v)^{q_1} + c_{34}\theta + d_\varrho \tag{5.13}$$

where $d_\varrho = c_{33}\Lambda_0\varrho^{q_\gamma^*} + c_{34}\Lambda_1\varrho^{q_1}$.

**Step 2:** Estimates related to $M_1, M_2$ and $M$.

From (2.9), we have

$$\begin{aligned}
M_1(\mathbb{S}_n(v)) &= \|\mathbb{G}_\mu[\delta^\gamma g_n(v+\varrho\mathbb{K}_\mu[\nu])]\|_{L_w^{q_\gamma^*}(\Omega;\delta^{\alpha_++\gamma})} \\
&\leq c_7\|g_n(v+\varrho\mathbb{K}_\mu[\nu])\|_{L^1(\Omega;\delta^{\alpha_++\gamma})}.
\end{aligned} \tag{5.14}$$

It follows that

$$M_1(\mathbb{S}_n(v)) \leq c_7c_{33}\Lambda_0 M_1(v)^{q_\gamma^*} + c_7c_{35}\Lambda_1 M_2(v)^{q_1} + c_7c_{34}\theta + c_7d_\varrho. \tag{5.15}$$

Applying (2.9), we get

$$\begin{aligned}
M_2(\mathbb{S}_n(v)) &= \|\mathbb{G}_\mu[\delta^\gamma g_n(v+\varrho\mathbb{K}_\mu[\nu])]\|_{L^{q_1}(\Omega;\delta^{\tilde{\gamma}})} \\
&\leq c_{36}\|g_n(v+\varrho\mathbb{K}_\mu[\nu])\|_{L^1(\Omega;\delta^{\alpha_++\gamma})},
\end{aligned}$$

which implies

$$M_2(\mathbb{S}_n(v)) \leq c_{36}c_{33}\Lambda_0 M_1(v)^{q_\gamma^*} + c_{36}c_{35}\Lambda_1 M_2(v)^{q_1} + c_{36}c_{34}\theta + c_{36}d_\varrho. \tag{5.16}$$

Consequently,

$$M(\mathbb{S}_n(v)) \leq c_{37}\Lambda_0 M_1(v)^{q_\gamma^*} + c_{38}\Lambda_1 M_2(v)^{q_1} + c_{40}\theta + c_{39}d_\varrho \tag{5.17}$$

where $c_{37} = (c_7 + c_{36})c_{33}$, $c_{38} = (c_7 + c_{36})c_{35}$, $c_{39} = c_7 + c_{36}$, $c_{40} = (c_7 + c_{36})c_{34}$. Therefore if $M(\mathbb{S}_n(v)) \le \lambda$ then

$$M(\mathbb{S}_n(v)) \le c_{37}\Lambda_0\lambda^{q^*_\gamma} + c_{38}\Lambda_1\lambda^{q_1} + c_{40}\theta + c_{39}d_\varrho.$$

Since $q^*_\gamma > 1$ and $q_1 > 1$, there exist $\varrho_0 > 0$ and $\theta_0 > 0$ such that for any $\varrho \in (0, \varrho_0)$ and $\theta \in (0, \theta_0)$ the equation

$$c_{37}\Lambda_0\lambda^{q^*_\gamma} + c_{38}\Lambda_1\lambda^{q_1} + c_{40}\theta + c_{39}d_\varrho = \lambda$$

admits a largest root $\bar\lambda > 0$. Therefore,

$$M(v) \le \bar\lambda \Longrightarrow M(\mathbb{S}_n(v)) \le \bar\lambda. \tag{5.18}$$

**Step 3:** We apply Schauder fixed point theorem to our setting.

Set

$$\mathcal{O} = \{\phi \in L^1_+(\Omega) : M(\phi) \le \bar\lambda\}.$$

Clearly, $\mathcal{O}$ is a convex subset of $L^1(\Omega)$. We shall show that $\mathcal{O}$ is a closed subset of $L^1(\Omega)$. Indeed, let $\{\phi_m\}$ be a sequence in $\mathcal{O}$ converging to $\phi$ in $L^1(\Omega)$. Obviously, $\phi \ge 0$. We can extract a subsequence, still denoted by $\{\phi_m\}$, such that $\phi_m \to \phi$ a.e. in $\Omega$. Consequently, by Fatou's lemma, $M_i(\phi) \le \liminf_{m\to\infty} M_i(\phi_m)$ for $i = 1, 2$. It follows that $M(\phi) \le \bar\lambda$. So $\phi \in \mathcal{O}$ and therefore $\mathcal{O}$ is a closed subset of $L^1(\Omega)$.

In light of (5.13) and (5.18), $\mathbb{S}_n$ is well-defined in $\mathcal{O}$ and $\mathbb{S}_n(\mathcal{O}) \subset \mathcal{O}$.

We observe that $\mathbb{S}_n$ is continuous. Indeed, if $\phi_m \to \phi$ as $m \to \infty$ in $L^1(\Omega)$ then $g_n(\phi_m + \varrho\mathbb{K}_\mu[\nu]) \to g_n(\phi + \varrho\mathbb{K}_\mu[\nu])$ as $m \to \infty$ in $L^1(\Omega; \delta^{\alpha_+ + \gamma})$. By (2.9), $\mathbb{S}_n(\phi_m) \to \mathbb{S}_n(\phi)$ as $m \to \infty$ in $L^1(\Omega)$.

We next show that $\mathbb{S}_n$ is a compact operator. Let $\{\phi_m\} \subset \mathcal{O}$ and for each $n$ put $\psi_m = \mathbb{S}_n(\phi_m)$. Hence $\{\Delta\psi_n\}$ is uniformly bounded in $L^p(G)$ for every compact subset $G \subset \Omega$. Therefore $\{\psi_m\}$ is uniformly bounded in $W^{1,p}(G)$. Consequently, there exists a subsequence, still denoted by $\{\psi_m\}$, and a function $\psi$ such that $\psi_m \to \psi$ a.e. in $\Omega$. By dominated convergence theorem, $\psi_m \to \psi$ in $L^1(\Omega)$. Thus $\mathbb{S}_n$ is compact.

By Schauder fixed point theorem there is a function $v_n \in L^1_+(\Omega)$ such that $\mathbb{S}_n(v_n) = v_n$ and $M(v_n) \le \bar\lambda$ where $\bar\lambda$ is independent of $n$. Due to Proposition 2.4, $\operatorname{tr}^*(v_n) = 0$ and $v_n$ is a nonnegative solution of (5.3). Moreover, there holds

$$-\int_\Omega v_n L_\mu \zeta dx = \int_\Omega \delta^\gamma g_n(v_n + \varrho\mathbb{K}_\mu[\nu])\zeta dx \quad \forall \zeta \in X(\Omega). \tag{5.19}$$

□

**Proof of Theorem G.** Let $\theta \in (0, \theta_0)$ and $\varrho \in (0, \varrho_0)$. For each $n$, set $u_n = v_n + \varrho\mathbb{K}_\mu[\nu]$ where $v_n$ is the solution constructed in Lemma 5.1. Then $\operatorname{tr}^*(u_n) = \varrho\nu$ and

$$-\int_\Omega u_n L_\mu \zeta dx = \int_\Omega \delta^\gamma g_n(u_n)\zeta dx - \varrho\int_\Omega \mathbb{K}_\mu[\nu] L_\mu\zeta dx \quad \forall \zeta \in X(\Omega). \tag{5.20}$$

Since $\{v_n\} \subset \mathcal{O}$, the sequence $\{g_n(v_n + \varrho\mathbb{K}_\mu[\nu])\}$ is uniformly bounded in $L^1(\Omega; \delta^{\alpha_+ + \gamma})$ and the sequence $\{\frac{\mu}{\delta^2} v_n\}$ is uniformly bounded in $L^{q_1}(G)$ for every compact subset $G \subset \Omega$. As a consequence, $\{\Delta v_n\}$ is uniformly bounded in $L^1(\Omega)$. By regularity result [18] for elliptic equations, there exists a subsequence, still denoted by $\{v_n\}$, and a function $v$ such that $v_n \to v$ a.e. in $\Omega$. Therefore $u_n \to u$ a.e. in $\Omega$ with $u = v + \varrho\mathbb{K}_\mu[\nu]$ and $g_n(u_n) \to \tilde g(u)$ a.e. in $\Omega$.

We show that $u_n \to u$ in $L^1(\Omega; \delta^{-\alpha_-})$. Since $\{v_n\}$ is uniformly bounded in $L^{q_1}(\Omega; \delta^{\tilde\gamma})$, by (2.10), we derive that $\{u_n\}$ is uniformly bounded in $L^{q_1}(\Omega; \delta^{-\alpha_-})$. Due to Holder inequality, $\{u_n\}$ is uniformly integrable with respect to $\delta^{-\alpha_-}dx$. We invoke Vitali's convergence theorem to derive that $u_n \to u$ in $L^1(\Omega; \delta^{-\alpha_-})$.

We next prove that $g_n(u_n) \to \tilde{g}(u)$ in $L^1(\Omega;\delta^{\alpha_++\gamma})$. For $\lambda > 0$ and $n \in \mathbb{N}$ set $B_{n,\lambda} = \{x \in \Omega : u_n > \lambda\}$ and $b_n(\lambda) = \int_{B_{n,\lambda}} \delta^{\alpha_++\gamma} dx$. For any Borel set $E \subset \Omega$,

$$
\begin{aligned}
\int_E g_n(u_n)\delta^{\alpha_++\gamma}dx &= \int_{E\cap B_{n,\lambda}} g_n(u_n)\delta^{\alpha_++\gamma}dx + \int_{E\cap B_{n,\lambda}^c} g_n(u_n)\delta^{\alpha_++\gamma}dx \\
&\le \int_{B_{n,\lambda}} g_n(u_n)\delta^{\alpha_++\gamma}dx + \Theta_\lambda \int_E \delta^{\alpha_++\gamma}dx \\
&\le b_n(\lambda)g_n(\lambda) + \int_\lambda^\infty b_n(s)dg_n(s) + \Theta_\lambda \int_E \delta^{\alpha_++\gamma}dx.
\end{aligned} \tag{5.21}
$$

where $\Theta_\lambda := \sup_{[0,\lambda]} g$. By proceeding as in the proof of Lemma 5.1, we deduce

$$
b_n(\lambda)g_n(\lambda) + \int_\lambda^\infty b_n(s)dg_n(s) \le c_{41}\int_\lambda^\infty s^{-1-q_\gamma^*} g_n(s)ds \le c_{41}\int_\lambda^\infty s^{-1-q_\gamma^*}\tilde{g}(s)ds \tag{5.22}
$$

where $c_{41}$ depends on $N,\mu,\gamma$ and $\Omega$. Note that the term on the right hand-side of (5.22) tends to 0 as $\lambda \to \infty$. Take arbitrarily $\varepsilon > 0$, there exists $\lambda > 0$ such that the right hand-side of (5.22) is smaller than $\frac{\varepsilon}{2}$. Fix such $\lambda$ and put $\eta = \frac{\varepsilon}{2\Theta_\lambda}$. Then, by (5.21),

$$
\int_E \delta(x)^{\alpha_++\gamma}dx \le \eta \Longrightarrow \int_E g_n(u_n)\delta(x)^{\alpha_++\gamma}dx < \varepsilon.
$$

Therefore the sequence $\{g_n(u_n)\}$ is uniformly integrable with respect to $\delta^{\alpha_++\gamma}dx$. Due to Vitali convergence theorem, we deduce that $g_n(u_n) \to \tilde{g}(u)$ in $L^1(\Omega;\delta^{\alpha_++\gamma})$.

Finally, by sending $n \to \infty$ in each term of (5.20) we obtain

$$
-\int_\Omega uL_\mu\zeta dx = \int_\Omega \delta^\gamma \tilde{g}(u)\zeta dx - \varrho\int_\Omega \mathbb{K}_\mu[\nu]L_\mu\zeta dx \quad \forall \zeta \in X(\Omega). \tag{5.23}
$$

By Theorem A, $u$ is a nonnegative weak solution of (1.26). □

5.2. **Sublinearity.** In this subsection we deal with the case where $g$ is sublinear.

**Lemma 5.2.** *Let Let $\nu \in \mathfrak{M}^+(\partial\Omega)$ such that $\|\nu\|_{\mathfrak{M}(\partial\Omega)} = 1$ and $\{g_n\} \subset C^1(\mathbb{R}_+)$ be a sequence satisfying* (5.1). *Assume* (1.29) *is satisfied. Then for every $\varrho > 0$ problem* (5.3) *admits a nonnegative solution $v_n$ satisfying*

$$
\|v_n\|_{L^1(\Omega;\delta^{\tilde{\gamma}})} \le \tilde{\lambda} \tag{5.24}
$$

*where $\tilde{\gamma}$ is as in* (5.2) *and $\tilde{\lambda}$ depends on $\Lambda_2, q_2, N, \mu$.*

**Proof.** The proof is similar to that of Lemma 5.1, also based on Schauder fixed point theorem. So we point out only the main modifications. Let $\mathbb{S}_n$ be the operator defined in (5.5). Fix $q_3 \in (1, \frac{N+\tilde{\gamma}}{N-1-\alpha_-})$ and put

$$
N(v) = \|v\|_{L^{q_3}(\Omega;\delta^{\tilde{\gamma}})}, \quad \forall v \in L^{q_3}(\Omega;\delta^{\tilde{\gamma}}).
$$

Combining (2.9), (2.10) and (1.29) leads to

$$\begin{aligned}N(\mathbb{S}_n(v)) &\leq c_{42}\left\|g_n(v+\varrho\mathbb{K}_\mu[\nu])\right\|_{L^1(\Omega;\delta^{\alpha_+ +\gamma})}\\ &\leq c_{42}\int_\Omega \Lambda_2(v+\varrho\mathbb{K}_\mu[\nu])^{q_2}\delta^{\alpha_+ +\gamma}dx + c_{42}\theta\int_\Omega \delta^{\alpha_+ +\gamma}dx\\ &\leq c_{42}\Lambda_2\int_\Omega v^{q_2}\delta^{\alpha_+ +\gamma}dx + c_{43}(\varrho^{q_2}+\theta)\\ &\leq c_{44}\Lambda_2\left(\int_\Omega v^{q_3}\delta^{\alpha_+ +\gamma}dx\right)^{\frac{q_2}{q_3}} + c_{43}(\varrho^{q_2}+\theta)\\ &\leq c_{44}\Lambda_2 N(v)^{q_2} + c_{43}(\varrho^{q_2}+\theta)\end{aligned}$$

where $c_i = c_i(N,\mu,\Omega,q_2)$ $(42 \leq i \leq 44)$. Therefore, if $N(v) \leq \lambda$ for some $\lambda > 0$ then

$$N(\mathbb{S}_n(v)) \leq c_{44}\Lambda_2\lambda^{q_2} + c_{43}(\varrho^{q_2}+\theta).$$

Consider the following algebraic equation

$$c_{44}\Lambda_2\lambda^{q_2} + c_{43}(\varrho^{q_2}+\theta) = \lambda. \tag{5.25}$$

If $q_2 < 1$ then for any $\varrho > 0$ (5.25) admits a unique positive root $\tilde{\lambda}$. If $q_2 = 1$ then for $\Lambda_2$ small such that $c_{44}\Lambda_2 < 1$ and $\varrho > 0$ equation (5.25) admits a unique positive root $\tilde{\lambda}$. Therefore,

$$N(v) \leq \tilde{\lambda} \Longrightarrow N(\mathbb{S}_n(v)) \leq \tilde{\lambda}. \tag{5.26}$$

By proceeding as in the proof of Lemma 5.1, we can prove that $\mathbb{S}_n$ is a continuous, compact operator from the closed, convex set

$$\tilde{\mathcal{O}} = \{v \in L^1_+(\Omega) : N(v) \leq \tilde{\lambda}\}$$

into itself. Thus by appealing to Schauder fixed point theorem, we see that there exists a function $v_n \in L^1_+(\Omega)$ such that $\mathbb{S}_n(v_n) = v_n$ and $N(v_n) \leq \tilde{\lambda}$ with $\tilde{\lambda}$ being independent of $n$. By Proposition 2.4, $\mathrm{tr}^*(v_n) = 0$ and $v_n$ is a nonnegative solution of (5.3). Moreover (5.19) holds. □

**Proof of Theorem H.** Let $v_n$ be the solution of (5.3) constructed in Lemma 5.2. Put $u_n = v_n + \varrho\mathbb{K}_\mu[\nu]$ then $u_n$ satisfies (5.20). By a similar argument as in the proof of Theorem G, there exists a subsequence, still denoted by $\{u_n\}$ and a function $u$ such that $u_n \to u$ a.e. in $\Omega$. Since $\{v_n\} \subset \tilde{\mathcal{O}}$, it follows that $\{v_n\}$ is uniformly bounded in $L^{q_3}(\Omega;\delta^{-\alpha_-})$, so is $\{u_n\}$. By Holder inequality, $\{u_n\}$ is uniformly integrable in $L^1(\Omega;\delta^{-\alpha_-})$. Due to (1.29), $\{g_n(u_n)\}$ is uniformly integrable in $L^1(\Omega;\delta^{\alpha_+ +\gamma})$. Vitali convergence theorem implies that $u_n \to u$ in $L^1(\Omega;\delta^{-\alpha_-})$ and $g_n(u_n) \to \tilde{g}(u)$ in $L^1(\Omega;\delta^{\alpha_+ +\gamma})$. Letting $n \to \infty$ in (5.20), we conclude that $u$ is a nonnegative solution of (1.9) satisfying (1.13). □

5.3. **Supercritical case. Proof of Theorem I.** Suppose by contradiction that for some $\varrho > 0$ and $z \in \partial\Omega$ there exists a positive solution $u$ of (1.31). Then by Theorem A, $\tilde{g}(u) \in L^q(\Omega;\delta^{\alpha_+ +\gamma})$ and

$$u = \mathbb{G}_\mu[\delta^\gamma \tilde{g}(u)] + \varrho K_\mu(\cdot, z).$$

This, along with (2.7) and the monotonicity assumption, implies

$$\begin{aligned}\int_\Omega \delta(x)^{\alpha_++\gamma}\tilde{g}(u)dx &\geq \int_\Omega \delta(x)^{\alpha_++\gamma}\tilde{g}(K_\mu(x,z))dx\\ &\geq \int_\Omega \delta(x)^{\alpha_++\gamma}\tilde{g}(\delta(x)^{\alpha_+}|x-z|^{2\alpha_--N})dx\\ &\geq c_{45}\int_{\{x\in\Omega:\delta(x)\geq\frac{1}{2}|x-z|\}}|x-z|^{\alpha_++\gamma}\tilde{g}(|x-z|^{1+\alpha_--N})dx.\end{aligned}$$

Fix $r_0>0$ such that

$$\mathcal{C}_{r_0}:=\left\{x:|x|\leq r_0,\delta(x)\geq\frac{1}{2}|x-z|\right\}\subset\left\{x\in\Omega:\delta(x)\geq\frac{1}{2}|x-z|\right\}.$$

Then

$$\begin{aligned}\int_\Omega \delta(x)^{\alpha_++\gamma}\tilde{g}(u)dx &\geq c_{45}\int_{\mathcal{C}_{r_0}}|x-z|^{\alpha_++\gamma}\tilde{g}(|x-z|^{1+\alpha_--N})dx\\ &\geq c_{46}\int_{a^*}^\infty s^{-1-q^*_\gamma}\tilde{g}(s)ds\end{aligned}$$

for some $a^*>0$. By assumption (1.30), $\tilde{g}(u)\notin L^1(\Omega;\delta^{\alpha_++\gamma})$, which leads to a contradiction. □